\documentclass[12pt]{article}

\usepackage{amsmath,amssymb,latexsym,theorem,bbm}

\newfont{\msa}{msam10 scaled\magstep1}
\newfont{\ssmsa}{msam9}
\newfont{\smsa}{msam10}
\newfont{\sms}{msbm10}
\newfont{\sseufb}{eufb9}
\newfont{\seufb}{eufb10}
\newfont{\eufb}{eufb10 scaled\magstep1}
\newfont{\eusb}{eusb10 scaled\magstep1}
\newfont{\hcmr}{cmr17 scaled\magstep5}

\newcommand{\AstkoneKn}{\raise-7mm\hbox{\hcmr*}%
          ^{\hspace*{-5.3mm}\raise2.8mm\hbox{$\scriptstyle K_n$}}%
           _{\hspace*{-6.1mm}\raise2.5mm\hbox{$\scriptstyle k=1$}}}
\newcommand{\Astkonen}{\raise-7mm\hbox{\hcmr*}%
          ^{\hspace*{-5.3mm}\raise2.8mm\hbox{$\scriptstyle n$}}%
           _{\hspace*{-6.1mm}\raise2.5mm\hbox{$\scriptstyle k=1$}}}

\newcommand{\DS}{\displaystyle}

\newcommand{\SC}{\scriptstyle}

\newcommand{\CC}{\mathbb{C}}

\newcommand{\EE}{\mathrm{E}}
\newcommand{\LL}{\mathbb{L}}

\newcommand{\NN}{\mathbb{N}}
\newcommand{\PP}{\mathrm{P}}

\newcommand{\RR}{\mathbb{R}}
\newcommand{\TT}{\mathbb{T}}

\newcommand{\ZZ}{\mathbb{Z}}

\newcommand{\sSS}{\raise-0.5truemm\hbox{\sms S}}

\newcommand{\cI}{\mathcal{I}}

\newcommand{\cC}{\mathcal{C}}
\newcommand{\cL}{\mathcal{L}}
\newcommand{\cP}{\mathcal{P}}

\newcommand{\cN}{\mathcal{N}}

\newcommand{\dd}{\mathrm{d}}
\newcommand{\ee}{\mathrm{e}}
\newcommand{\qq}{\mathrm{q}}
\newcommand{\uu}{\mathrm{u}}
\newcommand{\ww}{\mathrm{w}}
\newcommand{\RE}{\mathrm{Re}\,}
\newcommand{\IM}{\mathrm{Im}\,}

\newcommand{\Var}{\mathrm{Var}}
\newcommand{\arc}{\mathrm{arc}}

\newcommand{\hmu}{\widehat{\mu}}

\newcommand{\homega}{\widehat{\omega}}
\newcommand{\hgamma}{\widehat{\gamma}}
\newcommand{\hpi}{\widehat{\pi}}

\newcommand{\tg}{\widetilde{g}}
\newcommand{\hDelta}{\widehat{\Delta}}
\newcommand{\hS}{\widehat{S}}
\newcommand{\hTT}{\widehat{\TT}}

\newcommand{\hG}{\widehat{G}}

\newcommand{\sleq}{\mbox{\ssmsa\hspace*{0.1mm}\symbol{54}\hspace*{0.1mm}}}

\renewcommand{\leq}{\mbox{\msa\hspace*{0.9mm}\symbol{54}\hspace*{0.9mm}}}
\renewcommand{\geq}{\mbox{\msa\hspace*{0.9mm}\symbol{62}\hspace*{0.9mm}}}

\newcommand{\bone}{\mathbbm{1}}
\newcommand{\vare}{\varepsilon}

\newcommand{\proofend}{\hfill\mbox{$\Box$}}

\newcommand{\distr}{\stackrel{\csD}{\longrightarrow}}
\newcommand{\distre}{\stackrel{\csD}{=}}
\newcommand{\weak}{\stackrel{\ww}{\longrightarrow}}

\newcommand{\csD}{{\SC\mathcal{D}}}

\numberwithin{equation}{section}

\theoremstyle{change} \theorembodyfont{\em}
\newtheorem{Lem}{Lemma.}[section]
\newtheorem{Thm}[Lem]{Theorem.}

\theorembodyfont{\rm}
\newtheorem{Def}[Lem]{Definition.}
\newtheorem{Rem}[Lem]{Remark.}

\begin{document}

\begin{center}
  {\bfseries\Large Limit theorems  \\
         on locally compact Abelian groups}\\ [5mm]
 {\sc M\'aty\'as Barczy},%
\footnote{Faculty of Informatics, University of Debrecen, Pf.12,
 H--4010 Debrecen, Hungary.
 E--mail: \mbox{barczy@inf.unideb.hu}}
 {\sc Alexander Bendikov}%
\footnote{Mathematical Institute, Wroclaw University, pl.
Grundwaldzki 2/4, 50-384 Wroclaw, Poland.
 E--mail: bendikov@math.uni.wroc.pl}
 and {\sc Gyula Pap}%
\footnote{Faculty of Informatics, University of Debrecen, Pf.12,
  H--4010 Debrecen, Hungary.
  E--mail: papgy@inf.unideb.hu,
  Phone: +36\,52\,512\,900/2826, Fax: +36\,52\,416\,857}
\end{center}

\vspace*{-2mm}

{\small\textbf{Abstract.} We prove limit theorems for row sums of a
rowwise independent infinitesimal array of random variables with
values in a locally compact Abelian group. First we give a proof of
Gaiser's theorem \cite[Satz 1.3.6]{GAI94}, since it does not have an
easy access and it is not complete. This theorem gives sufficient
conditions for convergence of the row sums, but the limit measure
can not have a nondegenerate idempotent factor. Then we prove
necessary and sufficient conditions for convergence of the
 row sums, where the limit measure can be also a nondegenerate Haar measure on
 a compact subgroup. Finally, we investigate special cases: the torus group, the group of
 \ $p$--adic integers and the \ $p$--adic solenoid.}

{\small\textbf{2000 Mathematics Subject Classification:} 60B10, 60B15}

\section{Introduction}

Let \ $G$ \ be a locally compact Abelian \ $T_0$--topological group
having a countable basis of its topology. The main question of limit
problems on \ $G$ \ can be formulated as follows. Let \
$\{X_{n,k}:n\in\NN,\,k=1,\ldots,K_n\}$ \ be an array of rowwise
 independent random elements with values in \ $G$ \ satisfying the
 infinitesimality condition
 \[
  \lim_{n\to\infty}\,\max_{1\sleq k\sleq K_n}\PP(X_{n,k}\in G\setminus U)=0
 \]
 for all Borel neighbourhoods \ $U$ \ of the identity \ $e$ \ of \ $G$.
\ One searches for conditions on the array so that the convergence
in distribution
 \[
  \sum_{k=1}^{K_n}X_{n,k}\distr\mu\qquad\text{as \ $n\to\infty$}
 \]
 to a probability measure \ $\mu$ \ on \ $G$ \ holds.

Let \ $\cL(G)$ \ denote the set of all possible limits of row sums
of rowwise independent infinitesimal triangular arrays in \ $G$. \
The following problems arise:
 \begin{enumerate}
  \item[(P1)] How to parametrize the set \ $\cL(G)$, \ i.e., to give a
               bijection between \ $\cL(G)$ \ and an appropriate parameter set
               \ $\cP(G)$;
  \item[(P2)] How to associate suitable quantities \ $q_n$ \ to the rows
               \ $\{X_{n,k}:1\leq k\leq K_n\}$, \ $n\in\NN$ \ so that
               \[
                \sum_{k=1}^{K_n}X_{n,k}\distr\mu
                \quad\Longleftrightarrow\quad
                q_n\to q,
               \]
              where \ $q\in\cP(G)$ \ corresponds to the limiting distribution
               \ $\mu$, \ and the convergence \ $q_n\to q$ \ is meant in
               an appropriate sense.
 \end{enumerate}
The problem (P1) has been solved by Parthasarathy (see Chapter IV,
Corollary 7.1 in \cite{PAR} and Remark \ref{unique} in Section
\ref{par}). It turns out that any measure \ $\mu\in\cL(G)$ \ is
necessarily weakly infinitely divisible. Gaiser \cite{GAI94} gave a
partial solution to the problem (P2). His theorem (see Section
\ref{Gaiserclt}) gives only some sufficient conditions
 for the convergence \ $\sum_{k=1}^{K_n}X_{n,k}\distr\mu$, \ which does not
 include the case where \ $\mu$ \ has a nondegenerate idempotent factor, i.e.,
 a nondegenerate Haar measure on a compact subgroup of \ $G$.
\ For a survey of results on limit theorems on a general locally
compact Abelian topological group see Bingham \cite{BIN}.

In this paper we prove necessary and sufficient conditions for some
limit theorems to hold on general locally compact Abelian groups.
Our results complete the results of the paper \cite{GAI94}. In our
theorems the limit measure can be also a nondegenerate Haar measure
on a compact subgroup of \ $G$.

We also specify our results considering some classical groups such
as the torus group, the group of \ $p$--adic integers and the \
$p$--adic solenoid. Here we apply Gaiser's theorem as well. For
completeness, we present a proof of this theorem, since Gaiser's
 dissertation does not have an easy access and Gaiser's proof is not complete.
Concerning limit problems on totally disconnected Abelian
topological groups, like the group of \ $p$--adic integers, we
mention Tel\"oken \cite{TEL}.

\section{Parametrization of weakly infinitely divisible measures}
\label{par}

Let \ $\NN$ \ and \ $\ZZ_+$ \ denote the sets of positive and of
nonnegative integers, respectively. The expression ``a measure \
$\mu$ \ on \ $G$'' means a measure \ $\mu$ \ on
 the \ $\sigma$--algebra of Borel subsets of \ $G$.
\ The Dirac measure at a point \ $x\in G$ \ will be denoted by \
$\delta_x$.

\begin{Def}
A probability measure \ $\mu$ \ on \ $G$ \ is called
 \emph{weakly infinitely divisible} if for all \ $n\in\NN$ \ there exist a
 probability measure \ $\mu_n$ \ on \ $G$ \ and an element \ $x_n\in G$ \ such
 that \ $\mu=\mu_n^{*n}*\delta_{x_n}$, \ where \ $\mu_n^{*n}$ \ denotes the
 \ $n$--times convolution.
The collection of all weakly infinitely divisible measures on \ $G$
\ will be
 denoted by \ $\cI_{\ww}(G)$.
\end{Def}

According to Parthasarathy \cite[Chapter IV, Corollary 7.1]{PAR},
 \ $\cL(G)\subset\cI_\ww(G)$.
\ Now we recall the building blocks of weakly infinitely divisible
measures. The main tool for their description is the Fourier
transform. The character group of \ $G$ \ will be denoted by \
$\hG$. \ For every bounded measure \ $\mu$ \ on \ $G$, \ let \
$\hmu:\hG\to\CC$ \ be
 defined by
 \[
  \hmu(\chi):=\int_G\chi\,\dd\mu,
  \qquad \chi\in\hG.
 \]
This function \ $\hmu$ \ is called the \emph{Fourier transform} of \
$\mu$. \ The basic properties of the Fourier transformation can be
found, e.g., in
 Heyer \cite[Theorem 1.3.8, Theorem 1.4.2]{HEY}, in Hewitt and Ross
 \cite[Theorem 23.10]{HR} and in Parthasarathy
 \cite[Chapter IV, Theorem 3.3]{PAR}.

If \ $H$ \ is a compact subgroup of \ $G$ \ then \ $\omega_H$ \ will
denote the
 Haar measure on \ $H$ \ (considered as a measure on \ $G$\,) \ normalized by
 the requirement \ $\omega_H(H)=1$.
\ The normalized Haar measures of compact subgroups of \ $G$ \ are
the only
 idempotents in the semigroup of probability measures on \ $G$ \ (see, e.g.,
 Wendel \cite[Theorem 1]{WEN}).
For all \ $\chi\in\hG$,
 \begin{equation}\label{homega}
  \homega_H(\chi)
  =\begin{cases}
    1 & \text{if \ $\chi(x)=1$ \ for all \ $x\in H$,}\\
    0 & \text{otherwise,}
   \end{cases}
 \end{equation}
 i.e., \ $\homega_H=\bone_{H^\perp},$ \ where
 \[
  H^\perp:=\big\{\chi\in\hG:\text{$\chi(x)=1$ \ for all \ $x\in H$}\big\}
 \]
 is the annihilator of \ $H$.
\ Clearly \ $\omega_H\in\cI_\ww(G)$, \ since \
$\omega_H*\omega_H=\omega_H$. \ Sazonov and Tutubalin
\cite{SAZ-TUT66} proved that \ $\omega_H\in\cL(G)$.

Obviously \ $\delta_x\in\cI_\ww(G)$ \ for all \ $x\in G$, \ and one
can easily
 check that \ $\delta_x\in\cL(G)$ \ for all \ $x\in G^\arc$, \ where \ $G^\arc$
 \ denotes the arc--component of the identity \ $e$.

A \emph{quadratic form} on \ $\hG$ \ is a nonnegative continuous
function
 \ $\psi:\hG\to\RR_+$ \ such that
 \[
  \psi(\chi_1\chi_2)+\psi(\chi_1\chi_2^{-1})=2(\psi(\chi_1)+\psi(\chi_2))
  \qquad\text{for all \ $\chi_1,\chi_2\in\hG$.}
 \]
The set of all quadratic forms on \ $\hG$ \ will be denoted by \
$\qq_+(\hG)$. \ For any quadratic form \ $\psi\in\qq_+(\hG)$, \
there exists a unique
 probability measure \ $\gamma_\psi$ \ on \ $G$ \ determined by
 \[
  \hgamma_\psi(\chi)=\ee^{-\psi(\chi)/2}\qquad\text{for all \ $\chi\in\hG$},
 \]
 which is a symmetric Gauss measure (see, e.g., Theorem 5.2.8 in Heyer
 \cite{HEY}).
Obviously \ $\gamma_\psi\in\cL(G)$, \ since
 \ $\gamma_\psi=\gamma_{\psi/n}^{*n}$ \ for all \ $n\in\NN$ \ and
 \ $\gamma_{\psi/n}\weak\delta_e$ \ as \ $n\to\infty$.
\ (Here and in the sequel \ $\weak$ \ denotes weak convergence of
bounded
 measures on \ $G$.)

For a bounded measure \ $\eta$ \ on \ $G$, \ the
 \emph{compound Poisson measure} \ $\ee(\eta)$ \ is the probability measure
 on \ $G$ \ defined by
 \[
  \ee(\eta)
  :=\ee^{-\eta(G)}
    \left(\delta_e+\eta+\frac{\eta*\eta}{2!}+\frac{\eta*\eta*\eta}{3!}
          +\cdots\right).
 \]
The Fourier transform of a compound Poisson measure \ $\ee(\eta)$ \
is
 \begin{equation}\label{cP}
  (\ee(\eta))\:\widehat{}\:(\chi)
  =\exp\left\{\int_G(\chi(x)-1)\,\dd\eta(x)\right\},
   \qquad\chi\in\hG.
 \end{equation}
Clearly \ $\ee(\eta)\in\cL(G)$, \ since
 \ $\ee(\eta)=\big(\ee(\eta/n)\big)^{*n}$ \ for all \ $n\in\NN$ \ and
 \ $\ee(\eta/n)\weak\delta_e$ \ as \ $n\to\infty$.
\ In order to introduce generalized Poisson measures, we recall the
notions of
 a local inner product and a L\'evy measure. A Borel neighbourhood
  \ $U$ \ of \ $e$ \ is a Borel subset of \ $G$ \ for which there
  exists an open subset \ $\widetilde U$ \ such that
   \ $e\in\widetilde U\subset U.$ \ Let \ $\cN_e$ \ denote the
 collection of all Borel neighbourhoods \mbox{of \ $e$}.

\begin{Def}\label{lip}
A continuous function \ $g:G\times\hG\to\RR$ \ is called a
 \emph{local inner product} for \ $G$ \ if
 \begin{enumerate}
  \item[(i)] for every compact subset \ $C$ \ of \ $\hG$, \ there exists
         \ $U\in\cN_e$ \ such that
         \[\chi(x)=\ee^{ig(x,\chi)}\qquad
           \text{for all \ $x\in U$,\quad$\chi\in C$,}\]
  \item[(ii)] for all \ $x\in G$ \ and \ $\chi,\chi_1,\chi_2\in\hG$,
         \[g(x,\chi_1\chi_2)=g(x,\chi_1)+g(x,\chi_2),\qquad
           g(-x,\chi)=-g(x,\chi),\]
  \item[(iii)] for every compact subset \ $C$ \ of \ $\hG$,
         \[\sup\limits_{x\in G}\sup\limits_{\chi\in C}|g(x,\chi)|<\infty,
           \qquad
           \lim\limits_{x\to e}\sup\limits_{\chi\in C}|g(x,\chi)|=0.\]
 \end{enumerate}
\end{Def}

Parthasarathy \cite[Chapter IV, Lemma 5.3]{PAR} proved the existence
of a local
 inner product for an arbitrary locally compact Abelian \ $T_0$--topological
 group having a countable basis of its topology.

\begin{Def}
An extended real--valued measure \ $\eta$ \ on \ $G$ \ is said to be
a \emph{L\'evy measure} if \ $\eta(\{e\})=0$,
 \ $\eta(G\setminus U)<\infty$
 \ for all \ $U\in\cN_e$, \ and \ $\int_G(1-\RE\chi(x))\,\dd\eta(x)<\infty$
 \ for all \ $\chi\in\hG$.
\ The set of all L\'evy measures on \ $G$ \ will be denoted by \
$\LL(G)$.
\end{Def}

We note that for all \ $\chi\in\hG$ \ there exists \ $U\in\cN_e$ \
such that
 \begin{equation}\label{compare}
  \frac{1}{4}g(x,\chi)^2\leq 1-\RE\chi(x)\leq\frac{1}{2}g(x,\chi)^2,
  \qquad x\in U,
 \end{equation}
 thus the requirement \ $\int_G(1-\RE\chi(x))\,\dd\eta(x)<\infty$ \ can be
 replaced by \ $\int_G g(x,\chi)^2\,\dd\eta(x)<\infty$ \ for some (and then
 necessarily for any) local inner product \ $g$.

For a L\'evy measure \ $\eta\in\LL(G)$ \ and for a local inner
product \ $g$
 \ for \ $G$, \ the \emph{generalized Poisson measure} \ $\pi_{\eta,\,g}$
 \ is the probability measure on \ $G$ \ defined by
 $$
  \hpi_{\eta,\,g}(\chi)
  =\exp\left\{\int_G\big(\chi(x)-1-ig(x,\chi)\big)\,\dd\eta(x)\right\}\qquad
  \text{for all \ $\chi\in\hG$}
 $$
 (see, e.g., Chapter IV, Theorem 7.1 in Parthasarathy \cite{PAR}).
Obviously \ $\pi_{\eta,\,g}\in\cL(G)$, \ since
 \ $\pi_{\eta,\,g}=\pi_{\eta/n,\,g}^{*n}$ \ for all \ $n\in\NN$ \ and
 \ $\pi_{\eta/n,\,g}\weak\delta_e$ \ as \ $n\to\infty$.

\begin{Def}
For a bounded measure \ $\eta$ \ on \ $G$ \ and for a local inner
product \ $g$ \ for \ $G$, \ the \emph{local mean} of \ $\eta$ \
with respect to \ $g$ \ is the uniquely defined element \
$m_g(\eta)\in G$ \ given by
 \[
  \chi(m_g(\eta))=\exp\left\{i\int_G g(x,\chi)\,\dd\eta(x)\right\}\qquad
  \text{for all \ $\chi\in\hG$.}
 \]
\end{Def}
The existence and uniqueness of a local mean is guaranteed by
Pontryagin's duality theorem. If \ $\eta$ \ coincides with the
distribution \ $\PP_X$ \ of a random element
 \ $X$ \ in \ $G$, \ we will use the notation \ $m_g(X)$ \ instead of
 \ $m_g(\PP_X)$.
\ Remark that \ $\chi(m_g(X))=\ee^{i\,\EE\,g(X,\chi)}$ \ for all
 \ $\chi\in\hG$.

Note that for a bounded measure \ $\eta$ \ on \ $G$ \ with \
$\eta(\{e\})=0$
 \ we have \ $\eta\in\LL(G)$ \ and
 \ $\ee(\eta)=\pi_{\eta,\,g}*\delta_{m_g(\eta)}$.

Let \ $\cP(G)$ \ be the set of all quadruplets \ $(H,a,\psi,\eta)$,
\ where
 \ $H$ \ is a compact subgroup of \ $G$, \ $a\in G$, \ $\psi\in\qq_+(\hG)$
 \ and \ $\eta\in\LL(G)$.
\ Parthasarathy \cite[Chapter IV, Corollary 7.1]{PAR} proved the
following
 parametrization for weakly infinitely divisible measures on \ $G$.

\begin{Thm}[Parthasarathy]\label{LCA1}
Let \ $g$ \ be a fixed local inner product for \ $G$. \ If \
$\mu\in\cI_\ww(G)$ \ then there exists a quadruplet
 \ $(H,a,\psi,\eta)\in\cP(G)$ \ such that
 \begin{equation}\label{LH}
  \mu=\omega_H*\delta_a*\gamma_\psi*\pi_{\eta,\,g}.
 \end{equation}
Conversely, if \ $(H,a,\psi,\eta)\in\cP(G)$ \ then
 \ $\omega_H*\delta_a*\gamma_\psi*\pi_{\eta,\,g}\in\cI_\ww(G)$.
\end{Thm}

\begin{Rem}\label{unique}
In general, this parametrization is not one--to--one (see
Parthasarathy
 \cite[p.112, Remark 3]{PAR}), but the compact subgroup \ $H$ \ is uniquely
 determined in \eqref{LH} by \ $\mu$ \ (more precisely, \ $H$ \ is the
 annihilator of the open subgroup \ $\{\chi\in\hG:\hmu(\chi)\not=0\}$).
\ If \ $H=\{e\}$ \ then the quadratic form \ $\psi$ \ in \eqref{LH}
is also
 uniquely determined by \ $\mu$.
\ In order to obtain one--to--one parametrization one can take
equivalence
 classes of quadruplets related to the equivalence relation \ $\approx$
 \ defined by
 \[
  (H,a_1,\psi_1,\eta_1)\approx(H,a_2,\psi_2,\eta_2)
  \quad\Longleftrightarrow\quad
  \omega_H*\delta_{a_1}*\gamma_{\psi_1}*\pi_{\eta_1,\,g}
  =\omega_H*\delta_{a_2}*\gamma_{\psi_2}*\pi_{\eta_2,\,g}.
 \]
\end{Rem}

\section{Gaiser's limit theorem}
\label{Gaiserclt}

For a sequence \ $\{X_n:n\in\NN\}$ \ of random elements in \ $G$ \
and for a
 probability measure \ $\mu$ \ on \ $G$, \ notation \ $X_n\distr\mu$ \ means
 weak convergence \ $\PP_{X_n}\weak\mu$ \ of the distributions \ $\PP_{X_n}$
 \ of \ $X_n$, \ $n\in\NN$ \ towards \ $\mu$.
\ Let \ $\cC(G)$, \ $\cC_0(G)$ \ and \ $\cC_0^\uu(G)$ \ denote the
spaces of
 real--valued bounded continuous functions on \ $G$, \ the set of all functions
 in \ $\cC(G)$ \ vanishing in some \ $U\in\cN_e$, \ and the set of all
 uniformly continuous functions in \ $\cC_0(G)$, \ respectively.
Gaiser \cite[Satz 1.3.6]{GAI94} proved the following limit theorem.

\begin{Thm}[Gaiser]\label{THM:GAISER}
Let \ $g$ \ be a fixed local inner product for \ $G$. \ Let \
$\{X_{n,k}:n\in\NN,\,k=1,\ldots,K_n\}$ \ be a rowwise independent
 infinitesimal array of random elements in \ $G$.
\ Suppose that there exists a quadruplet \
$(\{e\},a,\psi,\eta)\in\cP(G)$
 \ such that
 \begin{enumerate}
 \item[(i)] \ $\DS\sum_{k=1}^{K_n}m_g(X_{n,k})\to a$ \ as \ $n\to\infty$,
 \item[(ii)] \ $\DS\sum_{k=1}^{K_n}\Var\,g(X_{n,k},\chi)
          \to\psi(\chi)+\int_Gg(x,\chi)^2\,\dd\eta(x)$
        \ as \ $n\to\infty$ \ for all \ $\chi\in\hG$,
 \item[(iii)] \ $\DS\sum_{k=1}^{K_n}\EE\,f(X_{n,k})\to\int_Gf\,\dd\eta$
        \ as \ $n\to\infty$ \ for all \ $f\in\cC_0(G)$.
 \end{enumerate}
Then
 \begin{equation}\label{Gaiser}
  \sum_{k=1}^{K_n}X_{n,k}\distr\delta_a*\gamma_\psi*\pi_{\eta,\,g}\qquad
  \text{as \ $n\to\infty$.}
 \end{equation}
\end{Thm}

\begin{Rem}
If either \ $a\not=e$ \ or \ $\psi\not=0$ \ or \ $\eta\not=0$ \ then
the infinitesimality of \ $\{X_{n,k}:n\in\NN,\,k=1,\ldots,K_n\}$ \
and \eqref{Gaiser} imply \ $K_n\to\infty$.
\end{Rem}

\begin{Rem}
Condition (i) is equivalent to
 \begin{enumerate}
  \item[$\text{(i}^\prime\text{)}$]
   $\DS\exp\left\{i\sum_{k=1}^{K_n}\EE\,g(X_{n,k},\chi)\right\}\to\chi(a)$
   \ as \ $n\to\infty$ \ for all \ $\chi\in\hG$.
 \end{enumerate}
\end{Rem}

Concerning condition (iii) we mention the following version of the
well-known portmanteau theorem (for the equivalence of (a) and (c),
see Meerschaert and Scheffler \cite[Proposition 1.2.19]{MS} and for
the rest, see Barczy and Pap \cite{BP}).

\begin{Thm}\label{Portmanteau}
Let \ $\{\eta_n:n\in\ZZ_+\}$ \ be a sequence of extended
real--valued measures
 on \ $G$ \ such that \ $\eta_n(G\setminus U)<\infty$ \ for all \ $U\in\cN_e$
 \ and for all \ $n\in\ZZ_+$.
\ Then the following assertions are equivalent:
 \begin{enumerate}
 \item[(a)] \ $\DS\int_Gf\,\dd\eta_n\to\int_Gf\,\dd\eta_0$ \ as \ $n\to\infty$
             \ for all \ $f\in\cC_0(G)$,
 \item[(b)] \ $\DS\int_Gf\,\dd\eta_n\to\int_Gf\,\dd\eta_0$ \ as \ $n\to\infty$
             \ for all \ $f\in\cC_0^\uu(G)$,
 \item[(c)] \ $\eta_n(G\setminus U)\to\eta_0(G\setminus U)$ \ as \ $n\to\infty$
        \ for all \ $U\in\cN_e$ \ with \ $\eta_0(\partial U)=0$,
 \item[(d)]
       \ $\DS\int_{G\setminus U}f\,\dd\eta_n\to\int_{G\setminus U}f\,\dd\eta_0$
        \ as \ $n\to\infty$ \ for all \ $f\in\cC(G)$, \ $U\in\cN_e$ \ with
        \ $\eta_0(\partial U)=0$,
 \item[(e)] \ $\eta_n|_{G\setminus U}\weak\eta_0|_{G\setminus U}$ \ as
        \ $n\to\infty$ \ for all \ $U\in\cN_e$ \ with \ $\eta_0(\partial U)=0$.
 \end{enumerate}
\end{Thm}

(Here and in the sequel \ $\eta|_B$ \ denotes the restriction of a
measure
 \ $\eta$ \ to a Borel subset \ $B$ \ of \ $G$, \ considered as a measure on
 \ $G$.)

\begin{Rem}\label{iii'}
Due to Theorem \ref{Portmanteau}, condition (iii) of Theorem
\ref{THM:GAISER}
 is equivalent to
 \begin{enumerate}
  \item[$\text{(iii}^\prime\text{)}$]
   $\DS\sum_{k=1}^{K_n}\PP(X_{n,k}\in G\setminus U)\to\eta(G\setminus U)$
   \ as \ $n\to\infty$ \ for all \ $U\in\cN_e$ \ with \ $\eta(\partial U)=0$.
 \end{enumerate}
\end{Rem}

In order to prove Theorem \ref{THM:GAISER}, first we recall a
theorem about
 convergence of weakly infinitely divisible measures without idempotent factors
 (see Gaiser \cite[Satz 1.2.1]{GAI94}).

\begin{Thm}\label{CWID}
For each \ $n\in\ZZ_+$, \ let \ $\mu_n\in\cI_\ww(G)$ \ be such that
\eqref{LH}
 holds for \ $\mu_n$ \ with a quadruplet \ $(\{e\},a_n,\psi_n,\eta_n)$.
\ If there exists a local inner product \ $g$ \ for \ $G$ \ such
that
 \begin{enumerate}
  \item[(i)] \ $a_n\to a_0$ \ as \ $n\to\infty$,
  \item[(ii)] \ $\psi_n(\chi)+\int_Gg(x,\chi)^2\,\dd\eta_n(x)
           \to\psi_0(\chi)+\int_Gg(x,\chi)^2\,\dd\eta_0(x)$
         \ as \ $n\to\infty$ \ for all \ $\chi\in\hG$,
  \item[(iii)] \ $\int_Gf\,\dd\eta_n\to\int_Gf\,\dd\eta_0$ \ as \ $n\to\infty$
         \ for all \ $f\in\cC_0(G)$,
 \end{enumerate}
 then \ $\mu_n\weak\mu_0$ \ as \ $n\to\infty$.
\end{Thm}

\noindent{\bf Proof.} It suffices to show \
$\hmu_n(\chi)\to\hmu_0(\chi)$ \ as \ $n\to\infty$ \ for
 all \ $\chi\in\hG$.
\ Let
 \[
  h(x,\chi):=\chi(x)-1-ig(x,\chi)+\frac{1}{2}g(x,\chi)^2
 \]
 for all \ $x\in G$ \  and all \ $\chi\in\hG$.
\ Then
 \[
  \hmu_n(\chi)
  =\chi(a_n)
   \exp\left\{-\frac{1}{2}
               \left(\psi_n(\chi)+\int_Gg(x,\chi)^2\,\dd\eta_n(x)\right)
              +\int_Gh(x,\chi)\,\dd\eta_n(x)\right\}
 \]
 for all \ $n\in\ZZ_+$ \ and all \ $\chi\in\hG$.
\ Taking into account assumptions (i) and (ii), it is enough to show
that
 \begin{equation}\label{h}
  \int_Gh(x,\chi)\,\dd\eta_n(x)\to\int_Gh(x,\chi)\,\dd\eta_0(x)\quad
  \text{as \ $n\to\infty$ \ for all \ $\chi\in\hG$.}
 \end{equation}
For each \ $\chi\in\hG$, \ there exists \ $U\in\cN_e$ \ such that
 \ $\chi(x)=\ee^{ig(x,\chi)}$ \ for all \ $x\in U$.
\ Using the inequality
 \begin{equation}\label{Taylor3}
  \left|\ee^{iy}-1-iy+\frac{y^2}{2}\right|\leq\frac{|y|^3}{6}\qquad
  \text{for all \ $y\in\RR$,}
 \end{equation}
 we obtain \ $|h(x,\chi)|\leq|g(x,\chi)|^3/6$ \ for all \ $x\in U$.
\ Consequently, for all \ $V\in\cN_e$ \ with \ $V\subset U$,
 \[
  \left|\int_Gh(x,\chi)\,\dd\eta_n(x)-\int_Gh(x,\chi)\,\dd\eta_0(x)\right|
  \leq I_n^{(1)}(V)+I_n^{(2)}(V),
 \]
 where
 \begin{align*}
  I_n^{(1)}(V)
  &:=\frac{1}{6}\int_V|g(x,\chi)|^3\,\dd(\eta_n+\eta_0)(x),\\[2mm]
  I_n^{(2)}(V)&:=\left|\int_{G\setminus V}h(x,\chi)\,\dd\eta_n(x)
                       -\int_{G\setminus V}h(x,\chi)\,\dd\eta_0(x)\right|.
 \end{align*}
We have
 \[
  I_n^{(1)}(V)
  \leq\frac{1}{6}\sup_{x\in V}|g(x,\chi)|
       \int_Vg(x,\chi)^2\,\dd(\eta_n+\eta_0)(x).
 \]
By assumption (ii),
 \[
  \sup_{n\in\ZZ_+}\int_Vg(x,\chi)^2\,\dd\eta_n(x)
  \leq\sup_{n\in\ZZ_+}
      \left(\psi_n(\chi)+\int_Gg(x,\chi)^2\,\dd\eta_n(x)\right)
  <\infty.
 \]
Theorem 8.3 in Hewitt and Ross \cite{HR} yields existence of a
metric \ $d$
 \ on \ $G$ \ compatible with the topology of \ $G$.
\ The function \ $t\mapsto\eta_0(\{x\in G:d(x,e)\geq t\})$ \ from
 \ $(0,\infty)$ \ into \ $\RR$ \ is non--increasing, hence the set
 \ $\big\{t\in(0,\infty):\eta_0(\{x\in G:d(x,e)=t\})>0\big\}$
 \ of its discontinuities is countable.
Consequently, for arbitrary \ $\vare>0$, \ there exists \ $t>0$ \
such that
 \ $V_1:=\{x\in G:d(x,e)<t\}\in\cN_e$, \ $V_1\subset U$,
 \ $\eta_0(\partial V_1)=0$ \ and
 \[
  \sup_{y\in V_1}|g(x,\chi)|
  <\frac{3\vare}{2\sup\limits_{n\in\ZZ_+}\int_Vg(x,\chi)^2\,\dd\eta_n(x)},
 \]
 thus \ $I_n^{(1)}(V_1)<\vare/2$.
\ By assumption (iii) and Theorem \ref{Portmanteau},
 \ $I_n^{(2)}(V_1)<\vare/2$ \ for all sufficiently large \ $n$, \ hence we
 obtain
 \[
  \left|\int_Gh(x,\chi)\,\dd\eta_n(x)
        -\int_Gh(x,\chi)\,\dd\eta_0(x)\right|
  <\vare
 \]
 for all sufficiently large \ $n$, \ which implies \eqref{h}.
\proofend

The notion of a special local inner product is also needed.

\begin{Def}
A local inner product \ $g$ \ for \ $G$ \ is called \emph{special}
if it is
 uniformly continuous in its first variable, i.e., if for all \ $\vare>0$
 \ there exists \ $U\in\cN_e$ \ such that \ $|g(x,\chi)-g(y,\chi)|<\vare$
 \ for all \ $x,y\in G$ \ with \ $x-y\in U$.
\end{Def}

Gaiser \cite[Satz 1.1.4]{GAI94} proved the existence of a special
local inner product for an arbitrary locally compact Abelian \
$T_0$--topological group having a countable basis of its topology.
The proof goes along the lines of the proof of the existence of a
local inner product in Heyer \cite[Theorem 5.1.10]{HEY}.

\noindent{\bf Proof.} \emph{(Proof of Theorem \ref{THM:GAISER})} \
First we show that it is enough to prove the statement for a special
local inner product, namely, if the statement is true for some local
inner product
 \ $g$, \ then it is true for any local inner product \ $\tg$.
\ Suppose that assumptions (i)--(iii) hold for \ $\tg$ \ with a
quadruplet
 \ $(\{e\},a,\psi,\eta)$.
\ We show that they hold for \ $g$ \ with the quadruplet
 \ $(\{e\},a+m_{g,\,\tg}(\eta),\psi,\eta)$, \ where the element
 \ $m_{g,\,\tg}(\eta)\in G$ \ is uniquely determined by
 \[
  \chi(m_{g,\,\tg}(\eta))
  =\exp\left\{i\int_G(g(x,\chi)-\tg(x,\chi))\,\dd\eta(x)\right\}\qquad
  \text{for all \ $\chi\in\hG$.}
 \]
(Note that \ $g(\cdot,\chi)-\tg(\cdot,\chi)\in\cC_0(G)$ \ can be
checked
 easily.)
Hence we want to prove
 \begin{enumerate}
 \item[$\text{(i}^\prime\text{)}$]
  \ $\DS\sum_{k=1}^{K_n}m_g(X_{n,k})\to a+m_{g,\,\tg}(\eta)$ \ as
     \ $n\to\infty$,
 \item[$\text{(ii}^\prime\text{)}$]
  \ $\DS\sum_{k=1}^{K_n}\Var\,g(X_{n,k},\chi)
        \to\psi(\chi)+\int_Gg(x,\chi)^2\,\dd\eta(x)$
  \ as \ $n\to\infty$ \ for all \ $\chi\in\hG$,
 \item[$\text{(iii}^\prime\text{)}$]
  \ $\DS\sum_{k=1}^{K_n}\EE\,f(X_{n,k})\to\int_Gf\,\dd\eta$
  \ as \ $n\to\infty$ \ for all \ $f\in\cC_0(G)$.
 \end{enumerate}

Clearly $\text{(iii}^\prime\text{)}$ holds, since it is identical
with
 assumption (iii).

By assumption (i), in order to prove $\text{(i}^\prime\text{)}$ we
have to show
 \[
  \chi\left(\sum_{k=1}^{K_n}m_g(X_{n,k})
            -\sum_{k=1}^{K_n}m_{\tg}(X_{n,k})\right)
  \to\chi(m_{g,\,\tg}(\eta))\qquad\text{for all \ $\chi\in\hG$.}
 \]
We have
 \begin{align*}
  &\chi\left(\sum_{k=1}^{K_n}m_g(X_{n,k})
             -\sum_{k=1}^{K_n}m_{\tg}(X_{n,k})\right)
   =\prod_{k=1}^{K_n}\frac{\chi(m_g(X_{n,k}))}{\chi(m_{\tg}(X_{n,k}))}
   =\prod_{k=1}^{K_n}\frac{\ee^{i\EE\,g(X_{n,k},\chi)}}
                          {\ee^{i\EE\,\tg(X_{n,k},\chi)}}\\[2mm]
  &=\exp\left\{i\sum_{k=1}^{K_n}
                \EE\big(g(X_{n,k},\chi)-\tg(X_{n,k},\chi)\big)\right\}
   \to\exp\left\{i\int_G(g(x,\chi)-\tg(x,\chi))\,\dd\eta(x)\right\},
 \end{align*}
 where we applied assumption (iii) with the function
 \ $g(\cdot,\chi)-\tg(\cdot,\chi)\in\cC_0(G)$.

By assumption (ii), in order to prove $\text{(ii}^\prime\text{)}$ we
have to
 show
 \begin{equation}\label{AB}
  \sum_{k=1}^{K_n}\Var\,g(X_{n,k},\chi)-\sum_{k=1}^{K_n}\Var\,\tg(X_{n,k},\chi)
  \to\int_G\big(g(x,\chi)^2-\tg(x,\chi)^2\big)\,\dd\eta(x)
 \end{equation}
 for all \ $\chi\in\hG$, \ where
 \ $g(\cdot,\chi)^2-\tg(\cdot,\chi)^2\in\cC_0(G)$ \ can be checked easily.
We have
 \[
  \sum_{k=1}^{K_n}\Var\,g(X_{n,k},\chi)-\sum_{k=1}^{K_n}\Var\,\tg(X_{n,k},\chi)
  =A_n-B_n,
 \]
 where
 \begin{align*}
  A_n&:=\sum_{k=1}^{K_n}\EE\,\big(g(X_{n,k},\chi)^2-\tg(X_{n,k},\chi)^2\big),\\
  B_n&:=\sum_{k=1}^{K_n}
         \big[(\EE\,g(X_{n,k},\chi))^2-(\EE\,\tg(X_{n,k},\chi))^2\big].
 \end{align*}
Applying assumption (iii) with the function
 \ $g(\cdot,\chi)^2-\tg(\cdot,\chi)^2\in\cC_0(G)$, \ we obtain
 \begin{equation}\label{An}
  A_n\to\int_G\big(g(x,\chi)^2-\tg(x,\chi)^2\big)\,\dd\eta(x).
 \end{equation}
Moreover,
 \[
  B_n=\sum_{k=1}^{K_n}
       \EE\big(g(X_{n,k},\chi)-\tg(X_{n,k},\chi)\big)\,
       \EE\big(g(X_{n,k},\chi)+\tg(X_{n,k},\chi)\big)
 \]
 implies
 \[
  |B_n|
  \leq\max_{1\sleq k\sleq K_n}
       \EE\big(|g(X_{n,k},\chi)|+|\tg(X_{n,k},\chi)|\big)\,
       \sum_{k=1}^{K_n}\EE\,|g(X_{n,k},\chi)-\tg(X_{n,k},\chi)|.
 \]
Using assumption (iii) with the function
 \ $|g(\cdot,\chi)-\tg(\cdot,\chi)|\in\cC_0(G)$, \ we get
 \begin{equation}\label{Bn1}
  \sum_{k=1}^{K_n}\EE\,|g(X_{n,k},\chi)-\tg(X_{n,k},\chi)|
  \to\int_G|g(x,\chi)-\tg(x,\chi)|\,\dd\eta(x).
 \end{equation}
Infinitesimality of \ $\{X_{n,k}:n\in\NN,\,k=1,\ldots,K_n\}$ \
implies
 \begin{equation}\label{infX}
  \max_{1\sleq k\sleq K_n}\EE\,|g(X_{n,k},\chi)|\to0\qquad
  \text{for all \ $\chi\in\hG$.}
 \end{equation}
Indeed,
 \[
  \max_{1\sleq k\sleq K_n}\EE\,|g(X_{n,k},\chi)|
  \leq\sup_{x\in U}|g(x,\chi)|
      +\sup_{x\in G}|g(x,\chi)|\cdot
       \max_{1\sleq k\sleq K_n}\PP(X_{n,k}\in G\setminus U)
 \]
 for all \ $U\in\cN_e$ \ and for all \ $\chi\in\hG$, \ and (iii) of Definition
 \ref{lip} implies \ $\sup_{x\in U}|g(x,\chi)|\to0$ \ as \ $U\to\{e\}$.
\ Clearly \eqref{Bn1} and \eqref{infX} imply \ $B_n\to0$, \ hence,
by
 \eqref{An}, we obtain \eqref{AB}.

We conclude that assumptions (i)--(iii) hold for
 \ $g$ \ with the quadruplet \ $(\{e\},a+m_{g,\,\tg}(\eta),\psi,\eta)$.
\ Since we supposed that the statement is true for \ $g$, \ we get
 \[
  \sum_{k=1}^{K_n}X_{n,k}
  \distr\delta_{a+m_{g,\,\tg}(\eta)}*\gamma_\psi*\pi_{\eta,\,g}.
 \]
Hence
 \begin{align*}
  \EE\chi\left(\sum_{k=1}^{K_n}X_{n,k}\right)
  &\to\chi(a+m_{g,\,\tg}(\eta))
      \exp\left\{-\frac{1}{2}\psi(\chi)
                 +\int_G\big(\chi(x)-1-ig(x,\chi)\big)\,\dd\eta(x)\right\}
   \\[2mm]
  &=\chi(a)
    \exp\left\{-\frac{1}{2}\psi(\chi)
                +\int_G\big(\chi(x)-1-i\tg(x,\chi)\big)\,\dd\eta(x)\right\}
 \end{align*}
 for all \ $\chi\in\hG$, \ which implies
 \[
  \sum_{k=1}^{K_n}X_{n,k}\distr\delta_a*\gamma_\psi*\pi_{\eta,\,\tg}.
 \]

Thus we may suppose that \ $g$ \ is a special local inner product.
Let \ $Y_{n,k}:=X_{n,k}-m_g(X_{n,k})$ \ for all \ $n\in\NN$, \
$k=1,\dots,K_n$. \ We show that \
$\{Y_{n,k}:n\in\NN,\,k=1,\ldots,K_n\}$ \ is an infinitesimal
 array of rowwise independent random elements in \ $G$, \ and
 \begin{enumerate}
 \item[$\text{(i}^{\prime\prime}\text{)}$]
  \ $\DS\sum_{k=1}^{K_n}m_g(Y_{n,k})\to e$ \ as \ $n\to\infty$,
 \item[$\text{(ii}^{\prime\prime}\text{)}$]
  \ $\DS\sum_{k=1}^{K_n}\EE\big(g(Y_{n,k},\chi)^2\big)
        \to\psi(\chi)+\int_Gg(x,\chi)^2\,\dd\eta(x)$
  \ as \ $n\to\infty$ \ for all \ $\chi\in\hG$,
 \item[$\text{(iii}^{\prime\prime}\text{)}$]
  \ $\DS\sum_{k=1}^{K_n}\EE\,f(Y_{n,k})\to\int_Gf\,\dd\eta$
  \ as \ $n\to\infty$ \ for all \ $f\in\cC_0(G)$.
 \end{enumerate}
Infinitesimality of \ $\{Y_{n,k}:n\in\NN,\,k=1,\ldots,K_n\}$ \ is
equivalent to
 \begin{equation}\label{infY}
  \max_{1\sleq k\sleq K_n}|\EE\,\chi(Y_{n,k})-1|\to0\qquad
  \text{for all \ $\chi\in\hG$.}
 \end{equation}
We have
 \begin{align*}
  |\EE\,\chi(Y_{n,k})-1|
  &=\left|\frac{\EE\,\chi(X_{n,k})}{\chi(m_g(X_{n,k}))}-1\right|
  =\left|\frac{\EE\,\chi(X_{n,k})}{\ee^{i\EE\,g(X_{n,k},\chi)}}-1\right|\\[2mm]
  &=\big|\EE\,\chi(X_{n,k})-\ee^{i\EE\,g(X_{n,k},\chi)}\big|
   \leq|\EE\,\chi(X_{n,k})-1|+\big|\ee^{i\EE\,g(X_{n,k},\chi)}-1\big|.
 \end{align*}
Infinitesimality of \ $\{X_{n,k}:n\in\NN,\,k=1,\ldots,K_n\}$ \
implies
 \begin{equation}\label{infX2}
  \max_{1\sleq k\sleq K_n}|\EE\,\chi(X_{n,k})-1|\to0\qquad
  \text{for all \ $\chi\in\hG$.}
 \end{equation}
Infinitesimality of \ $\{X_{n,k}:n\in\NN,\,k=1,\ldots,K_n\}$ \
implies \eqref{infX} as well, hence using the inequality
 \ $|\ee^{iy}-1|\leq|y|$ \ for all \ $y\in\RR$, \ we get
 \[
  \max_{1\sleq k\sleq K_n}\big|\ee^{i\EE\,g(X_{n,k},\chi)}-1\big|\to0\qquad
  \text{for all \ $\chi\in\hG$,}
 \]
 and we obtain \eqref{infY}.

For $\text{(i}^{\prime\prime}\text{)}$, it is enough to show
 \[
  \sum_{k=1}^{K_n}\EE\,g(Y_{n,k},\chi)\to0\qquad\text{for all \ $\chi\in\hG$.}
 \]
Let \ $\chi\in\hG$ \ be fixed. Infinitesimality of \
$\{X_{n,k}:n\in\NN,\,k=1,\ldots,K_n\}$ \ implies that for
 all \ $V\in\cN_e$ \ and for all sufficiently large \ $n$ \ we have
 \ $m_g(X_{n,k})\in V$ \ for \ $k=1,\dots,K_n$.
\ Consequently, using \eqref{infX} as well, we conclude that for all
 sufficiently large \ $n$ \ we have
 \begin{equation}\label{lmX}
  g(m_g(X_{n,k}),\chi)=\EE\,g(X_{n,k},\chi)\qquad\text{for \ $k=1,\dots,K_n$.}
 \end{equation}
Infinitesimality of \ $\{X_{n,k}:n\in\NN,\,k=1,\ldots,K_n\}$ \ and
properties
 of the local inner product \ $g$ \ imply also the existence of \ $U\in\cN_e$
 \ such that \ $\eta(\partial U)=0$ \ and
 \begin{equation}\label{nbhd}
  g(x-m_g(X_{n,k}),\chi)-g(x,\chi)
  =-g(m_g(X_{n,k}),\chi)\qquad\text{for \ $x\in U$, \ $k=1,\dots,K_n$}
 \end{equation}
 for all sufficiently large \ $n$ \ (see Parthasarathy \cite[page 91]{PAR}).
Consequently, for all sufficiently large \ $n$, \ we obtain
 \begin{align*}
  \left|\sum_{k=1}^{K_n}\EE\,g(Y_{n,k},\chi)\right|
  &=\left|\sum_{k=1}^{K_n}
           \EE\,\Big(g(Y_{n,k},\chi)-g(X_{n,k},\chi)+g(m_g(X_{n,k}),\chi)\Big)
                \bone_{G\setminus U}(X_{n,k})\right|\\
  &\leq\Big(\max_{1\sleq k\sleq K_n}
             \sup_{x\in G}|g(x-m_g(X_{n,k}),\chi)-g(x,\chi)|\Big)\,
       \sum_{k=1}^{K_n}\PP(X_{n,k}\in G\setminus U)\\
  &\phantom{\leq}
       +\max_{1\sleq k\sleq K_n}|g(m_g(X_{n,k}),\chi)|\,
        \sum_{k=1}^{K_n}\PP(X_{n,k}\in G\setminus U)\to0.
 \end{align*}
Indeed,
 \begin{equation}\label{uc}
  \max_{1\sleq k\sleq K_n}\sup_{x\in G}|g(x-m_g(X_{n,k}),\chi)-g(x,\chi)|
  \to0\qquad\text{as \ $n\to\infty$,}
 \end{equation}
 since \ $g$ \ is  uniformly continuous in its first variable and for all
 \ $V\in\cN_e$ \ and for all sufficiently large \ $n$ \ we have
 \ $m_g(X_{n,k})\in V$ \ for \ $k=1,\dots,K_n$.
\ Moreover, \eqref{infX} and \eqref{lmX} imply
 \begin{equation}\label{inflc}
  \max\limits_{1\sleq k\sleq K_n}|g(m_g(X_{n,k}),\chi)|\to0\qquad
  \text{as \ $n\to\infty$},
 \end{equation}
 and assumption (iii) implies
 \begin{equation}\label{iiim}
  \sup\limits_{n\in\NN}\sum\limits_{k=1}^{K_n}\PP(X_{n,k}\in G\setminus U)
  <\infty.
 \end{equation}
To prove $\text{(ii}^{\prime\prime}\text{)}$, we have to show
 \[
  \sum_{k=1}^{K_n}
   \Big(\EE\,\big(g(Y_{n,k},\chi)^2\big)-\Var\,g(X_{n,k},\chi)\Big)
  \to0\qquad\text{for all \ $\chi\in\hG$.}
 \]
Consider again a neighbourhood \ $U\in\cN_e$ \ such that \
$\eta(\partial U)=0$
 \ and \eqref{nbhd} holds for all sufficiently large \ $n$.
\ We have
 \[
  \EE\,\big(g(Y_{n,k},\chi)^2\big)-\Var\,g(X_{n,k},\chi)
  =C_{n,k}+D_{n,k},
 \]
 where
 \begin{align*}
  C_{n,k}&:=\EE\,\Big(g(Y_{n,k},\chi)^2-g(X_{n,k},\chi)^2\Big)\bone_U(X_{n,k})
            +\big(\EE\,g(X_{n,k},\chi)\big)^2,\\
  D_{n,k}&:=\EE\,\Big(g(Y_{n,k},\chi)^2-g(X_{n,k},\chi)^2\Big)
                 \bone_{G\setminus U}(X_{n,k}).
 \end{align*}
For all sufficiently large \ $n$ \ we have \eqref{lmX}, hence
 \begin{align*}
  C_{n,k}&=\EE\,\Big(\big(g(X_{n,k},\chi)-g(m_g(X_{n,k}),\chi)\big)^2
                     -g(X_{n,k},\chi)^2\Big)\bone_U(X_{n,k})
           +\big(\EE\,g(X_{n,k},\chi)\big)^2\\
         &=g(m_g(X_{n,k}),\chi)^2\,\PP(X_{n,k}\in U)
           -2g(m_g(X_{n,k}),\chi)\,
            \EE\,\big(g(X_{n,k},\chi)\bone_U(X_{n,k})\big)\\
         &\phantom{=}
           +\big(\EE\,g(X_{n,k},\chi)\big)^2\\
         &=2\EE\,g(X_{n,k},\chi)
           \,\EE\,\big(g(X_{n,k},\chi)\bone_{G\setminus U}(X_{n,k})\big)
           -\big(\EE\,g(X_{n,k},\chi)\big)^2\,\PP(X_{n,k}\in G\setminus U).
 \end{align*}
Consequently, again by \eqref{lmX},
 \begin{equation}\label{Cn}
  |C_{n,k}|\leq\PP(X_{n,k}\in G\setminus U)
               \left(2|\EE\,g(X_{n,k},\chi)|\,\sup_{x\in G}|g(x,\chi)|
                     +|\EE\,g(X_{n,k},\chi)|^2\right).
 \end{equation}
Moreover,
 \[
  D_{n,k}=\EE\,\big(g(Y_{n,k},\chi)-g(X_{n,k},\chi)\big)
               \big(g(Y_{n,k},\chi)+g(X_{n,k},\chi)\big)
               \bone_{G\setminus U}(X_{n,k}),
 \]
 thus
 \begin{equation}\label{Dn}
  |D_{n,k}|\leq2\PP(X_{n,k}\in G\setminus U)\sup_{x\in G}|g(x,\chi)|
               \max_{1\sleq k\sleq K_n}\,
               \sup_{x\in G}\big|g(x-m_g(X_{n,k}),\chi)-g(x,\chi)\big|.
 \end{equation}
Now \eqref{Cn} and \eqref{Dn}, using \eqref{uc}, \eqref{inflc} and
 \eqref{iiim}, imply $\text{(ii}^{\prime\prime}\text{)}$.

To prove $\text{(iii}^{\prime\prime}\text{)}$, it is enough to show
 \begin{equation}\label{iiiprime}
  \sum_{k=1}^{K_n}\EE\,f(Y_{n,k})-\sum_{k=1}^{K_n}\EE\,f(X_{n,k})\to0
  \end{equation}
 for all \ $f\in\cC_0^\uu(G)$ \ (see Theorem \ref{Portmanteau}).
Choose \ $V\in\cN_e$ \ such that \ $f(x)=0$ \ for all \ $x\in V$. \
Then choose \ $U\in\cN_e$ \ such that \ $U-U\subset V$, \ where
 \ $U-U:=\{x-y:x,y\in U\}$.
\ Infinitesimality of \ $\{X_{n,k}:n\in\NN,\,k=1,\ldots,K_n\}$ \
implies that
 for all sufficiently large \ $n$ \ we have \ $m_g(X_{n,k})\in U$ \ for
 \ $k=1,\dots,K_n$, \ hence
 \ $f(Y_{n,k})-f(X_{n,k})
     =\big(f(Y_{n,k})-f(X_{n,k})\big)\bone_{G\setminus U}(X_{n,k}).$
 \ Consequently,
 \[
  \left|\sum_{k=1}^{K_n}\EE\,f(Y_{n,k})-\sum_{k=1}^{K_n}\EE\,f(X_{n,k})\right|
  \leq\sup_{x\in G}\big|f(x-m_g(X_{n,k}))-f(x)\big|\,
                   \sum_{k=1}^{K_n}\PP(X_{n,k}\in G\setminus U),
 \]
 and uniform continuity of \ $f$ \ and \eqref{iiim} imply \eqref{iiiprime}.

Now consider the shifted compound Poisson measures
 \[
  \nu_n
  :=\ee\bigg(\sum_{k=1}^{K_n}\PP_{Y_{n,k}}\bigg)
    *\delta_{\sum_{k=1}^{K_n}m_g(X_{n,k})},\qquad n\in\NN.
 \]
Clearly \ $\nu_n\in\cI_\ww(G)$ \ such that \eqref{LH} holds for \
$\nu_n$
 \ with the quadruplet
 \[
  \left(\{e\},\,\sum_{k=1}^{K_n}m_g(X_{n,k})+\sum_{k=1}^{K_n}m_g(Y_{n,k}),\,
        0,\,\sum_{k=1}^{K_n}\PP_{Y_{n,k}}\right).
 \]
By Theorem \ref{CWID}, using (i) and
 $\text{(i}^{\prime\prime}\text{)}$--$\text{(iii}^{\prime\prime}\text{)}$, we
 obtain
 \[
  \nu_n\weak\delta_a*\gamma_\psi*\pi_{\eta,\,g}.
 \]
Applying a theorem about the accompanying Poisson array due to
Parthasarathy
 \cite[Chapter IV, Theorem 5.1]{PAR}, we conclude the statement.
\proofend

\section{Limit theorems for symmetric arrays}
\label{Haarclt1}

The following theorem is an easy consequence of Theorem
\ref{THM:GAISER}.

\begin{Thm}[CLT for symmetric array]\label{THM:GenSymmetric1}
Let \ $g$ \ be a fixed local inner product for \ $G$. \ Let \
$\{X_{n,k}:n\in\NN,\,k=1,\ldots,K_n\}$ \ be a rowwise independent
array
 of random elements in \ $G$ \ such that \ $X_{n,k}\distre-X_{n,k}$ \ for all
 \ $n\in\NN$, \ $k=1,\ldots,K_n$.
\ Suppose that there exists a quadratic form \ $\psi$ \ on \ $\hG$ \
such that
 \begin{enumerate}
 \item[(i)] \ $\DS\sum_{k=1}^{K_n}\Var\,g(X_{n,k},\chi)\to\psi(\chi)$
        \ as \ $n\to\infty$ \ for all \ $\chi\in\hG$,
 \item[(ii)] \ $\DS\sum_{k=1}^{K_n}\PP(X_{n,k}\in G\setminus U)\to0$
        \ as \ $n\to\infty$ \ for all \ $U\in\cN_e$.
 \end{enumerate}
Then the array \ $\{X_{n,k}:n\in\NN,\,k=1,\ldots,K_n\}$ \ is
infinitesimal and
 \[
  \sum_{k=1}^{K_n}X_{n,k}\distr\gamma_\psi\qquad\text{as \ $n\to\infty$.}
 \]
\end{Thm}

The next theorem gives necessary and sufficient conditions in case
of a rowwise
 independent and identically distributed (i.i.d.) symmetric array.
It turns out that in this special case conditions of Theorem
 \ref{THM:GenSymmetric1} are not only sufficient but necessary as well.
If \ $G$ is compact then the limit measure can be the normalized
Haar measure
 on \ $G$.

\begin{Thm}[Limit theorem for rowwise i.i.d.\ symmetric array]
\label{THM:GenSymmetric2} Let \
$\{X_{n,k}:n\in\NN,\,k=1,\ldots,K_n\}$ \ be an infinitesimal,
rowwise
 i.i.d.\ array of random elements in \ $G$ \ such that \ $K_n\to\infty$
 \ and \ $X_{n,k}\distre-X_{n,k}$ \ for all \ $n\in\NN$, \ $k=1,\ldots,K_n$.

If \ $g$ \ is a local inner product for \ $G$ \ and \ $\psi$ \ is a
quadratic
 form on \ $\hG$, \ then the following statements are equivalent:
 \begin{enumerate}
  \item[(i)] $\DS\sum_{k=1}^{K_n}X_{n,k}\distr\gamma_\psi$ \ as \ $n\to\infty$,
  \item[(ii)] $K_n\big(1-\RE\EE\,\chi(X_{n,1})\big)\to\frac{\psi(\chi)}{2}$
         \ as \ $n\to\infty$ \ for all \ $\chi\in\hG$,
  \item[(iii)] $K_n\,\Var\,g(X_{n,1},\chi)\to\psi(\chi)$ \ as \ $n\to\infty$
         \ for all \ $\chi\in\hG$ \ and
         \ $K_n\,\PP(X_{n,1}\in G\setminus U)\to0$ \ as
         \ $n\to\infty$ \ for all \ $U\in\cN_e$.
 \end{enumerate}

If \ $G$ \ is compact then
 \[
  \sum_{k=1}^{K_n}X_{n,k}\distr\omega_G
  \quad\Longleftrightarrow\quad
  \text{$K_n\big(1-\RE\EE\,\chi(X_{n,1})\big)\to\infty$ \ for all
         \ $\chi\in\hG\setminus\{\bone_G\}$.}
 \]
\end{Thm}

For the proof of Theorem \ref{THM:GenSymmetric2}, we need the
following simple
 observation.

\begin{Lem}\label{exp}
Let \ $\{\alpha_n:n\in\NN\}$ \ be a sequence of real numbers such
that
 \ $\alpha_n\geq-n$ \ for all sufficiently large \ $n$, \ and let
 \ $\alpha\in\RR\cup\{-\infty,\infty\}$.
\ Then
 \[
  \left(1+\frac{\alpha_n}{n}\right)^n\to\ee^\alpha
  \qquad\Longleftrightarrow\qquad
  \alpha_n\to\alpha,
 \]
 where \ $\ee^{-\infty}:=0$ \ and \ $\ee^\infty:=\infty$.
\end{Lem}

\noindent{\bf Proof.} If \ $\alpha_n\to\alpha\in\RR$ \ then \
$\alpha_n/n\to0$, \ hence L'Hospital's
 rule gives
 \[
  \log\left[\left(1+\frac{\alpha_n}{n}\right)^n\right]
  =\alpha_n\cdot\frac{\log\left(1+\alpha_n/n\right)}{\alpha_n/n}
  \to\alpha.
 \]
If \ $\alpha_n\to-\infty$ \ then we choose \ $n_0\in\NN$ \ such that
 \ $\alpha_n\geq-n$ \ for all \ $n\geq n_0$, \ hence \ $1+\alpha_n/n\geq0$
 \ for all \ $n\geq n_0$, \ implying
 \ $\liminf\limits_{n\to\infty}(1+\alpha_n/n)^n\geq0$.
\ For each \ $M\in\RR$ \ there exists \ $n_M\in\NN$ \ such that
 \ $\alpha_n\leq M$ \ for all \ $n\geq n_M$.
\ Then \ $(1+\alpha_n/n)^n\leq(1+M/n)^n$ \ for all \
$n\geq\max(n_0,n_M)$,
 \ which implies
 \[
  \limsup_{n\to\infty}\left(1+\frac{\alpha_n}{n}\right)^n
  \leq\limsup_{n\to\infty}\left(1+\frac{M}{n}\right)^n
  =\ee^M.
 \]
Since \ $M$ \ is arbitrary, we obtain
 \ $\limsup\limits_{n\to\infty}(1+\alpha_n/n)^n\leq0$, \ and finally
 \ $\lim\limits_{n\to\infty}(1+\alpha_n/n)^n=0$.
\ The case of \ $\alpha_n\to\infty$ \ can be handled similarly.

If \ $(1+\alpha_n/n)^n\to\ee^\alpha$ \ and \ $\alpha_n\not\to\alpha$
\ then
 there exist subsequences \ $(n^\prime)$ \ and \ $(n^{\prime\prime})$ \ and
 \ $\alpha^\prime,\alpha^{\prime\prime}\in\RR\cup\{-\infty,\infty\}$ \ with
 \ $\alpha^\prime\not=\alpha^{\prime\prime}$ \ such that
 \ $\alpha_{n^\prime}\to\alpha^\prime$ \ and
 \ $\alpha_{n^{\prime\prime}}\to\alpha^{\prime\prime}$.
Then \
$(1+\alpha_{n^\prime}/n^\prime)^{n^\prime}\to\ee^{\alpha^\prime}$ \
and
 \ $(1+\alpha_{n^{\prime\prime}}/n^{\prime\prime})^{n^{\prime\prime}}
    \to\ee^{\alpha^{\prime\prime}}$
 \ lead to a contradiction.
\proofend

\noindent{\bf Proof.} \emph{(Proof of Theorem
\ref{THM:GenSymmetric2})} \ (i) $\Longleftrightarrow$ (ii).
Statement (i) is equivalent to
 \begin{equation}\label{(i)}
  \EE\,\chi\bigg(\sum_{k=1}^{K_n}X_{n,k}\bigg)\to\hgamma_\psi(\chi)\qquad
  \text{for all \ $\chi\in\hG$.}
 \end{equation}
We have \ $\hgamma_\psi(\chi)=\ee^{-\psi(\chi)/2}$. \ Clearly \
$X_{n,k}\distre-X_{n,k}$ \ implies
 \ $\EE\,\chi(X_{n,k})=\RE\EE\,\chi(X_{n,k})$, \ hence
 \begin{equation}\label{sym}
  \EE\,\chi\bigg(\sum_{k=1}^{K_n}X_{n,k}\bigg)
  =\big(\RE\EE\,\chi(X_{n,1})\big)^{K_n}
  =\left(1+\frac{K_n\big(\RE\EE\,\chi(X_{n,1})-1\big)}{K_n}\right)^{K_n}.
 \end{equation}
Infinitesimality of \ $\{X_{n,k}:n\in\NN,\,k=1,\ldots,K_n\}$ \
implies
 \ $\EE\,\chi(X_{n,1})\to1$ \ (see \eqref{infY}), thus
 \ $\RE\EE\,\chi(X_{n,1})-1\geq-1$ \ for all sufficiently large \ $n\in\NN$.
\ Hence by \ $K_n\to\infty$ \ and by Lemma \ref{exp} we conclude
that
 \eqref{(i)} and (ii) are equivalent.

(ii) $\Longrightarrow$ (iii). We have already proved that (ii)
implies (i), hence, by Theorem 5.4.2 in
 Heyer \cite{HEY}, (ii) implies \ $K_n\,\PP(X_{n,1}\in G\setminus U)\to0$ \ for
 all \ $U\in\cN_e$.
\ Clearly \ $X_{n,k}\distre-X_{n,k}$ \ implies \
$\EE\,g(X_{n,k},\chi)=0$,
 \ thus \ $\Var\,g(X_{n,1},\chi)=\EE\,\big(g(X_{n,1},\chi)^2\big)$.
\ Consequently, it is enough to show
 \begin{equation}\label{(iii)-(ii)}
  K_n\left(\RE\EE\,\chi(X_{n,1})-1
           +\frac{1}{2}\EE\,\big(g(X_{n,1},\chi)^2\big)\right)
  \to0\qquad\text{for all \ $\chi\in\hG$.}
 \end{equation}
For \ $\chi\in\hG$, \ choose \ $U\in\cN_e$ \ such that
 \ $\chi(x)=\ee^{ig(x,\chi)}$ \ and \eqref{compare} hold for all \ $x\in U$.
\ Then
 \[
  K_n\left(\RE\EE\,\chi(X_{n,1})-1
           +\frac{1}{2}\EE\,\big(g(X_{n,1},\chi)^2\big)\right)
  =A_n+B_n,
 \]
 where
 \begin{align*}
  A_n&:=K_n\,\RE\EE\,\left(\ee^{ig(X_{n,1},\chi)}-1-ig(X_{n,1},\chi)
                           +\frac{1}{2}g(X_{n,1},\chi)^2\right)
                     \bone_U(X_{n,1}),\\
  B_n&:=K_n\,\RE\EE\,\left(\chi(X_{n,1})-1
                           +\frac{1}{2}g(X_{n,1},\chi)^2\right)
                     \bone_{G\setminus U}(X_{n,1}).
 \end{align*}
By formulas \eqref{Taylor3} and \eqref{compare} we get
 \[
  |A_n|\leq\frac{1}{6}K_n\,\EE\,\big(|g(X_{n,1},\chi)|^3\,\bone_U(X_{n,1})\big)
  \leq\frac{4\big(K_n(1-\RE\EE\,\chi(X_{n,1}))\big)^{3/2}}{3K_n^{1/2}},
 \]
 hence \ $K_n\to\infty$ \ and assumption (ii) yield \ $A_n\to0$.
\ Moreover,
 \[
 |B_n|\leq\left(2+\frac{1}{2}\sup_{x\in G}g(x,\chi)^2\right)
          K_n\,\PP(X_{n,1}\in G\setminus U)
     \to0,
 \]
 thus we obtain \eqref{(iii)-(ii)}.

(iii) $\Longrightarrow$ (i) follows from Theorem
\ref{THM:GenSymmetric1}.

Convergence \ $\sum_{k=1}^{K_n}X_{n,k}\distr\omega_G$ \ is
equivalent to
 \begin{equation}\label{last}
  \EE\,\chi\bigg(\sum_{k=1}^{K_n}X_{n,k}\bigg)\to\homega_G(\chi)\qquad
  \text{for all \ $\chi\in\hG$.}
 \end{equation}
Using \eqref{sym}, \eqref{homega} and Lemma \ref{exp}, one can
easily show that
 \eqref{last} holds if and only if
 \ $K_n\big(1-\RE\EE\,\chi(X_{n,1})\big)\to\infty$ \ for all
 \ $\chi\in\hG\setminus\{\bone_G\}$.
\proofend

The next statement is a special case of Theorem
\ref{THM:GenSymmetric2}.

\begin{Thm}[Limit theorem for rowwise i.i.d.\ Rademacher array]
\label{THM:GenRademacher} Let \ $x_n\in G$, \ $n\in\NN$ \ such that
\ $x_n\to e$. \  Let \ $\{X_{n,k}:n\in\NN,\,k=1,\ldots,K_n\}$ \ be a
rowwise i.i.d.\ array of
 random elements in \ $G$ \ such that \ $K_n\to\infty$ \ and
 \[
  \PP(X_{n,k}=x_n)=\PP(X_{n,k}=-x_n)=\frac{1}{2}.
 \]
Then the array \ $\{X_{n,k}:n\in\NN,\,k=1,\ldots,K_n\}$ \ is
infinitesimal.

If \ $\psi$ \ is a quadratic form on \ $\hG$ \ then
 \[
  \sum_{k=1}^{K_n}X_{n,k}\distr\gamma_\psi
  \quad\Longleftrightarrow\quad
  \text{$K_n\big(1-\RE\chi(x_n)\big)\to\frac{\psi(\chi)}{2}$ \ for all
         \ $\chi\in\hG$.}
 \]
\indent If \ $G$ \ is compact then
 \[
  \sum_{k=1}^{K_n}X_{n,k}\distr\omega_G
  \quad\Longleftrightarrow\quad
  \text{$K_n\big(1-\RE\chi(x_n)\big)\to\infty$ \ for all
         \ $\chi\in\hG\setminus\{\bone_G\}$.}
 \]
\end{Thm}

Note that in Theorem \ref{THM:GenRademacher} the expression \
$1-\RE\chi(x_n)$
 \ can be replaced in both places by \ $\frac{1}{2}g(x_n,\chi)^2$, \ where
 \ $g$ \ is an arbitrary local inner product for \ $G$ \ (see the proof of
 \eqref{(iii)-(ii)} and the inequalities in \eqref{compare}).

\section{Limit theorem for Bernoulli arrays}
\label{Haarclt2}

In the following limit theorem the limit measure can be the
normalized Haar
 measure on an arbitrary compact subgroup of \ $G$ \ generated by a single
 element.

\begin{Thm}[Limit theorem for rowwise i.i.d.\ Bernoulli array]
\label{THM:GenBernoulli} Let \ $x\in G$ \ such that \ $x\not=e$. \
Let \ $\{X_{n,k}:n\in\NN,\,k=1,\ldots,K_n\}$ \ be a rowwise i.i.d.\
array of
 random elements in \ $G$ \ such that \ $K_n\to\infty$,
 \[
  \PP(X_{n,k}=x)=p_n,\qquad
  \PP(X_{n,k}=e)=1-p_n,
 \]
 and \ $p_n\to0$.
\ Then the array \ $\{X_{n,k}:n\in\NN,\,k=1,\ldots,K_n\}$ \ is
infinitesimal.

If \ $\lambda$ \ is a nonnegative real number then
 \[
  \sum_{k=1}^{K_n}X_{n,k}\distr\ee(\lambda\delta_x)
  \quad\Longleftrightarrow\quad
  K_n\,p_n\to\lambda.
 \]
\indent If the smallest closed subgroup \ $H$ \ of \ $G$ \
containing \ $x$ \ is
 compact then
 \[
  \sum_{k=1}^{K_n}X_{n,k}\distr\omega_H
  \quad\Longleftrightarrow\quad
  K_n\,p_n\to\infty.
 \]
\end{Thm}

\noindent{\bf Proof.} First we suppose \ $K_n\,p_n\to\lambda$ \ and
show convergence
 \ $\sum_{k=1}^{K_n}X_{n,k}\distr\ee(\lambda\delta_x)$.
\ We need to prove
 \begin{equation}\label{Poisson}
  \EE\,\chi\bigg(\sum_{k=1}^{K_n}X_{n,k}\bigg)\to
  (\ee(\lambda\delta_x))\:\widehat{}\:(\chi)\qquad
  \text{for all \ $\chi\in\hG$.}
 \end{equation}
We have \
$(\ee(\lambda\delta_x))\:\widehat{}\:(\chi)=\ee^{\lambda(\chi(x)-1)}$
 \ and
 \begin{equation}\label{symP}
  \EE\,\chi\bigg(\sum_{k=1}^{K_n}X_{n,k}\bigg)
  =(p_n\chi(x)+1-p_n)^{K_n}
  =\left(1+\frac{K_n\,p_n(\chi(x)-1)}{K_n}\right)^{K_n}.
 \end{equation}
If \ $\{z_n:n\in\NN\}$ \ is a sequence of complex numbers such that
 \ $z_n\to z\in\CC$ \ then \ $\left(1+\frac{z_n}{n}\right)^n\to\ee^z$.
\ Consequently, \ $K_n\,p_n\to\lambda$ \ and \ $K_n\to\infty$ \
imply
 \eqref{Poisson}.

Next we suppose \ $K_n\,p_n\to\infty$ \ and show that
 \ $\sum_{k=1}^{K_n}X_{n,k}\distr\omega_H$.
\ We need to prove
 \[
  \EE\,\chi\bigg(\sum_{k=1}^{K_n}X_{n,k}\bigg)\to\homega_H(\chi)\qquad
  \text{for all \ $\chi\in\hG$.}
 \]
According to Hewitt and Ross \cite[Remarks 23.24]{HR}, \
$\{x\}^\perp=H^\perp$,
 \ and thus by \eqref{homega} we are left to check
 \begin{equation}\label{omegaH}
  \EE\,\chi\bigg(\sum_{k=1}^{K_n}X_{n,k}\bigg)
  \to\begin{cases}
      1 & \text{if \ $\chi\in\{x\}^\perp$,}\\
      0 & \text{otherwise.}
     \end{cases}
 \end{equation}
If \ $\chi\in\{x\}^\perp$ \ then \ $\chi(x)=1$, \ hence
 \[
  \EE\,\chi\bigg(\sum_{k=1}^{K_n}X_{n,k}\bigg)=(p_n\chi(x)+1-p_n)^{K_n}=1,
 \]
 and we obtain \eqref{omegaH}.
To handle the case \ $\chi\not\in\{x\}^\perp$, \ consider the
equality
 \begin{align*}
  \left|\EE\,\chi\bigg(\sum_{k=1}^{K_n}X_{n,k}\bigg)\right|
  &=|p_n\chi(x)+1-p_n|^{K_n}
  =\Big(\big(1+p_n(\RE\chi(x)-1)\big)^2
        +p_n^2\big(\IM\chi(x)\big)^2\Big)^{K_n/2}\\
  &=\left(1+\frac{K_n\,p_n\Big(2(\RE\chi(x)-1)+p_n|1-\chi(x)|^2\Big)}
                 {K_n}\right)^{K_n/2}.
 \end{align*}
Clearly \ $\chi\not\in\{x\}^\perp$ \ implies \ $\chi(x)\not=1$, \
and by
 \ $|\chi(x)|=1$ \ we get \ $\RE\chi(x)-1<0$.
\ Hence, by Lemma \ref{exp}, we conclude that \ $K_n\,p_n\to\infty$,
 \ $K_n\to\infty$ \ and \ $p_n\to0$ \ imply \eqref{omegaH}.

Now we suppose \ $\sum_{k=1}^{K_n}X_{n,k}\distr\ee(\lambda\delta_x)$
\ and
 derive \ $K_n\,p_n\to\lambda$.
\ If \ $K_n\,p_n\not\to\lambda$ \ then either there exists a
subsequence
 \ $(n^\prime)$ \ such that \ $K_{n^\prime}\,p_{n^\prime}\to\infty$, \ or there
 exist subsequences \ $(n^{\prime\prime})$ \ and \ $(n^{\prime\prime\prime})$
 \ and two distinct nonnegative real numbers \ $\lambda^{\prime\prime}$ \ and
 \ $\lambda^{\prime\prime\prime}$ \ such that
 \ $K_{n^{\prime\prime}}\,p_{n^{\prime\prime}}\to\lambda^{\prime\prime}$ \ and
 \ $K_{n^{\prime\prime\prime}}\,p_{n^{\prime\prime\prime}}
    \to\lambda^{\prime\prime\prime}$.
\ In the first case we would obtain
 \ $\sum_{k=1}^{K_{n^\prime}}X_{n^\prime,k}\distr\omega_H$, \ which
 contradicts to \ $\sum_{k=1}^{K_n}X_{n,k}\distr\ee(\lambda\delta_x)$.
\ In the second case we would obtain
 \ $\sum_{k=1}^{K_{n^{\prime\prime}}}X_{n^{\prime\prime},k}
    \distr\ee(\lambda^{\prime\prime}\delta_x)$
 \ and
 \ $\sum_{k=1}^{K_{n^{\prime\prime\prime}}}X_{n^{\prime\prime\prime},k}
    \distr\ee(\lambda^{\prime\prime\prime}\delta_x)$
 \ which again contradicts to
 \ $\sum_{k=1}^{K_n}X_{n,k}\distr\ee(\lambda\delta_x)$.

Finally we suppose \ $\sum_{k=1}^{K_n}X_{n,k}\distr\omega_H$ \ and
prove
 \ $K_n\,p_n\to\infty$.
\ If \ $K_n\,p_n\not\to\infty$ \ then there exists a subsequence \
$(n^\prime)$
 \ and a nonnegative real number \ $\lambda^\prime$ \ such that
 \ $K_{n^\prime}\,p_{n^\prime}\to\lambda^\prime$.
\ Then we would obtain
 \ $\sum_{k=1}^{K_{n^\prime}}X_{n^\prime,k}\distr\ee(\lambda^\prime\delta_x)$,
 \ which contradicts to \ $\sum_{k=1}^{K_n}X_{n,k}\distr\omega_H$.
\proofend

\section{Limit theorems on the torus}
\label{gencltT}

The set \ $\TT:=\{\ee^{ix}:-\pi\leq x<\pi\}$ \ equipped with the
usual
 multiplication of complex numbers is a compact Abelian \ $T_0$--topological
 group having a countable basis of its topology.
This is called the 1--dimensional torus group. Its character group
is \ $\hTT=\{\chi_\ell:\ell\in\ZZ\}$, \ where
 \[
  \chi_\ell(y):=y^\ell,\qquad y\in\TT,\quad\ell\in\ZZ.
 \]
Hence \ $\hTT$ \ can be identified with the additive group of
integers \ $\ZZ$. \ The compact subgroups of \ $\TT$ \ are
 \[
  H_r:=\{\ee^{2\pi ij/r}:j=0,1,\dots,r-1\},\qquad r\in\NN,
 \]
 and \ $\TT$ \ itself.

The set of all quadratic forms on \ $\hTT\cong\ZZ$ \ is
 \ $\qq_+\big(\hTT\big)=\{\psi_b:b\in\RR_+\}$, \ where
 \[
  \psi_b(\chi_\ell):=b\ell^2,\qquad\ell\in\ZZ,\quad b\in\RR_+.
 \]
Let us define the functions \ $\arg:\TT\to[-\pi,\pi[$ \ and \
$h:\RR\to\RR$
 \ by
 \begin{align*}
  \arg(\ee^{ix})&:=x,\qquad-\pi\leq x<\pi,\\
  h(x)&:=\begin{cases}
               0 & \text{if \ $x<-\pi$ \ or \ $x\geq\pi$,}\\
          -x-\pi & \text{if \ $-\pi\leq x<-\pi/2$,}\\
               x & \text{if \ $-\pi/2\leq x<\pi/2$,}\\
          -x+\pi & \text{if \ $\pi/2\leq x<\pi$.}\\
         \end{cases}
 \end{align*}
An extended real--valued measure \ $\eta$ \ on \ $\TT$ \ is a L\'evy
measure if
 and only if \ $\eta(\{e\})=0$ \ and \ $\int_\TT(\arg y)^2\,\dd\eta(y)<\infty$.
\ The function \ $g_\TT:\TT\times\hTT\to\RR$, \ defined as
 \[
  g_\TT(y,\chi_\ell):=\ell h(\arg y),\qquad y\in\TT,\quad\ell\in\ZZ,
 \]
 is a local inner product for \ $\TT$.

Theorem \ref{THM:GAISER} has the following consequence on the torus.

\begin{Thm}[Gauss--Poisson limit theorem]\label{THM:GPT}
Let \ $\{X_{n,k}:n\in\NN,\,k=1,\ldots,K_n\}$ \ be a rowwise
independent array
 of random elements in \ $\TT$.
\ Suppose that there exists a quadruplet \
$(\{e\},a,\psi_b,\eta)\in\cP(\TT)$
 \ such that
 \begin{enumerate}
 \item[(i)] \ $\DS\max_{1\sleq k\sleq K_n}\PP(|\arg(X_{n,k})|>\vare)\to0$
        \ as \ $n\to\infty$ \ for all \ $\vare>0$,
 \item[(ii)]
       \ $\DS\exp\left\{i\sum_{k=1}^{K_n}\EE\,h(\arg(X_{n,k}))\right\}\to a$
        \ as \ $n\to\infty$,
 \item[(iii)] \ $\DS\sum_{k=1}^{K_n}\Var\,h(\arg(X_{n,k}))
          \to b+\int_\TT\big(h(\arg y)\big)^2\,\dd\eta(y)$
        \ as \ $n\to\infty$,
 \item[(iv)] \ $\DS\sum_{k=1}^{K_n}\EE\,f(X_{n,k})\to\int_\TT f\,\dd\eta$
        \ as \ $n\to\infty$ \ for all \ $f\in\cC_0(\TT)$.
 \end{enumerate}
Then the array \ $\{X_{n,k}:n\in\NN,\,k=1,\ldots,K_n\}$ \ is
infinitesimal and
 \[
  \sum_{k=1}^{K_n}X_{n,k}\distr\delta_a*\gamma_{\psi_b}*\pi_{\eta,\,g_\TT}
  \qquad\text{as \ $n\to\infty$.}
 \]
\end{Thm}

If the limit measure has no generalized Poisson factor \
$\pi_{\eta,\,g_\TT}$
 \ then the truncation function \ $h$ \ can be omitted.

\begin{Thm}[CLT]\label{THM:GT}
Let \ $\{X_{n,k}:n\in\NN,\,k=1,\ldots,K_n\}$ \ be a rowwise
independent array
 of random elements in \ $\TT$.
\ Suppose that there exist an element \ $a\in\TT$ \ and a
nonnegative real
 number \ $b$ \ such that
 \begin{enumerate}
 \item[(i)] \ $\DS\exp\left\{i\sum_{k=1}^{K_n}\EE\,\arg(X_{n,k})\right\}\to a$
        \ as \ $n\to\infty$,
 \item[(ii)] \ $\DS\sum_{k=1}^{K_n}\Var\,\arg(X_{n,k})\to b$ \ as
        \ $n\to\infty$,
 \item[(iii)] \ $\DS\sum_{k=1}^{K_n}\PP(|\arg(X_{n,k})|>\vare)\to0$ \ as
        \ $n\to\infty$ \ for all \ $\vare>0$.
 \end{enumerate}
Then the array \ $\{X_{n,k}:n\in\NN,\,k=1,\ldots,K_n\}$ \ is
infinitesimal and
 \[
  \sum_{k=1}^{K_n}X_{n,k}\distr\delta_a*\gamma_{\psi_b}\qquad
  \text{as \ $n\to\infty$.}
 \]
\end{Thm}

\noindent{\bf Proof.} In view of Theorem \ref{THM:GPT} and Remark
\ref{iii'}, it is enough to check
 \begin{enumerate}
 \item[$\text{(i}^\prime\text{)}$]
  \ $\DS\exp\left\{i\sum_{k=1}^{K_n}\EE\,h(\arg(X_{n,k}))\right\}\to a$ \ as
   \ $n\to\infty$,
 \item[$\text{(ii}^\prime\text{)}$]
  \ $\DS\sum_{k=1}^{K_n}\Var\,h(\arg(X_{n,k}))\to b$ \ as \ $n\to\infty$,
 \item[$\text{(iii}^\prime\text{)}$]
  \ $\DS\sum_{k=1}^{K_n}\PP(|\arg(X_{n,k})|>\vare)\to0$ \ as
        \ $n\to\infty$ \ for all \ $\vare>0$.
 \end{enumerate}
Clearly $\text{(iii}^\prime\text{)}$ and assumption (iii) are
identical. In order to prove $\text{(i}^\prime\text{)}$ it is
sufficient to show
 \[
  \sum_{k=1}^{K_n}\EE\,h(\arg(X_{n,k}))-\sum_{k=1}^{K_n}\EE\,\arg(X_{n,k})\to0,
 \]
 since \ $|\ee^{iy_1}-\ee^{iy_2}|=|\ee^{i(y_1-y_2)}-1|\leq|y_1-y_2|$ \ for all
 \ $y_1,y_2\in\RR$.
\ We have \ $|h(y)-y|\leq\pi\bone_{[-\pi,-\pi/2]\cup[\pi/2,\pi]}(y)$
\ for all
 \ $y\in[-\pi,\pi]$, \ hence
 \[
  \left|\sum_{k=1}^{K_n}\EE\,h(\arg(X_{n,k}))
        -\sum_{k=1}^{K_n}\EE\,\arg(X_{n,k})\right|
  \leq\pi\sum_{k=1}^{K_n}\PP(|\arg(X_{n,k})|\geq\pi/2)\to0
 \]
 by condition (iii).
In order to check $\text{(ii}^\prime\text{)}$ it is enough to prove
 \[
  \sum_{k=1}^{K_n}\Var\,h(\arg(X_{n,k}))-\sum_{k=1}^{K_n}\Var\,\arg(X_{n,k})
  \to0.
 \]
We have
 \begin{align*}
  &\left|\sum_{k=1}^{K_n}\Var\,h(\arg(X_{n,k}))
        -\sum_{k=1}^{K_n}\Var\,\arg(X_{n,k})\right|\\
  &\leq\sum_{k=1}^{K_n}\EE\,\Big|\big(h(\arg(X_{n,k}))\big)^2
                                 -(\arg(X_{n,k}))^2\Big|
       +\sum_{k=1}^{K_n}
         \Big|\big(\EE\,h(\arg(X_{n,k}))\big)^2-(\EE\,\arg(X_{n,k}))^2\Big|\\
  &\leq2\pi^2\sum_{k=1}^{K_n}\PP(|\arg(X_{n,k})|\geq\pi/2)\to0,
 \end{align*}
 as desired.
\proofend

Theorem \ref{THM:GenRademacher} has the following consequence on the
torus.

\begin{Thm}[Limit theorem for rowwise i.i.d.\ Rademacher array]
\label{THM:TRademacherT} Let \ $x_n\in\TT$, \ $n\in\NN$ \ such that
\ $x_n\to e$. \  Let \ $\{X_{n,k}:n\in\NN,\,k=1,\ldots,K_n\}$ \ be a
rowwise i.i.d.\ array of
 random elements in \ $\TT$ \ such that \ $K_n\to\infty$ \ and
 \[
  \PP(X_{n,k}=x_n)=\PP(X_{n,k}=-x_n)=\frac{1}{2}.
 \]
Then the array \ $\{X_{n,k}:n\in\NN,\,k=1,\ldots,K_n\}$ \ is
infinitesimal.

If \ $b$ \ is a nonnegative real number then
 \[
  \sum_{k=1}^{K_n}X_{n,k}\distr\gamma_{\psi_b}
  \qquad\Longleftrightarrow\qquad
  K_n(\arg x_n)^2\to b.
 \]
\indent Moreover,
 \[
  \sum_{k=1}^{K_n}X_{n,k}\distr\omega_\TT
  \qquad\Longleftrightarrow\qquad
  K_n(\arg x_n)^2\to\infty.
 \]
\end{Thm}

\section{Limit theorems on the group of \ $p$--adic integers}
\label{gencltD}

Let \ $p$ \ be a prime. The group of \ $p$--adic integers is
 \[
  \Delta_p
  :=\big\{(x_0,x_1,\dots)
          :\text{$x_j\in\{0,1,\dots,p-1\}$ \ for all \ $j\in\ZZ_+$}\big\},
 \]
 where the sum \ $z:=x+y\in\Delta_p$ \ for \ $x,y\in\Delta_p$ \ is uniquely
 determined by the relationships
 \[
  \sum_{j=0}^dz_jp^j\equiv\sum_{j=0}^d(x_j+y_j)p^j\quad\mod p^{d+1}
  \qquad\text{for all \ $d\in\ZZ_+$.}
 \]
For each \ $r\in\ZZ_+$, \ let
 \ $\Lambda_r:=\{x\in\Delta_p:\text{$x_j=0$ \ for all
     \ $j\leq r-1$}\}.$ \ The family of sets
  \ $\{x+\Lambda_r:x\in\Delta_p,\,r\in\ZZ_+\}$ \
 is an open subbasis for a topology on \ $\Delta_p$ \ under which
 \ $\Delta_p$ \ is a compact, totally disconnected Abelian \ $T_0$--topological group
having a countable basis of its topology. Its character group is
 \ $\hDelta_p=\{\chi_{d,\ell}:d\in\ZZ_+,\,\ell=0,1,\dots,p^{d+1}-1\}$, \ where
 \[
  \chi_{d,\ell}(x):=\ee^{2\pi i\ell(x_0+px_1+\cdots+p^dx_d)/p^{d+1}},\qquad
  x\in\Delta_p,\quad d\in\ZZ_+,\quad\ell=0,1,\dots,p^{d+1}-1.
 \]
The compact subgroups of \ $\Delta_p$ \ are \ $\Lambda_r$, \
$r\in\ZZ_+$
 \ (see Hewitt and Ross \cite[Example 10.16 (a)]{HR}).

An extended real--valued measure \ $\eta$ \ on \ $\Delta_p$ \ is a
L\'evy
 measure if and only if \ $\eta(\{e\})=0$ \ and
 \ $\eta(\Delta_p\setminus\Lambda_r)<\infty$ \ for all \ $r\in\ZZ_+$.

Since the group \ $\Delta_p$ \ is totally disconnected, the only
quadratic
 form on \ $\hDelta_p$ \ is \ $\psi=0$, \ and the function
 \ $g_{\Delta_p}:\Delta_p\times\hDelta_p\to\RR$, \ $g_{\Delta_p}=0$ \ is a
 local inner product for \ $\Delta_p$.

Theorem \ref{THM:GAISER} has the following consequence on the group
 \ $\Delta_p$ \ of \ $p$--adic integers.

\begin{Thm}[Poisson limit theorem]\label{THM:PD}
Let \ $\{X_{n,k}:n\in\NN,\,k=1,\ldots,K_n\}$ \ be a rowwise
independent array of random elements in \ $\Delta_p$. \ Suppose that
there exists a L\'evy measure \ $\eta\in\LL(\Delta_p)$ \ such
 that
 \begin{enumerate}
 \item[(i)] $\DS\max_{1\sleq k\sleq K_n}
        \PP\Big(\big((X_{n,k})_0,\dots,(X_{n,k})_d\big)\not=0\Big)
        \to0$
       \ as \ $n\to\infty$ \ for all \ $d\in\ZZ_+$,
 \item[(ii)] $\DS\sum_{k=1}^{K_n}
        \PP\big((X_{n,k})_0=\ell_0,\dots,(X_{n,k})_d=\ell_d\big)
        \to\eta(\{x\in\Delta_p:x_0=\ell_0,\dots,x_d=\ell_d\})$
      \ as \ $n\to\infty$ \ for all \ $d\in\ZZ_+$,
      \ $\ell_0,\dots,\ell_d\in\{0,\dots,p-1\}$ \ with
      \ $(\ell_0,\dots,\ell_d)\not=0$.
 \end{enumerate}
Then the array \ $\{X_{n,k}:n\in\NN,\,k=1,\ldots,K_n\}$ \ is
infinitesimal and
 \[
  \sum_{k=1}^{K_n}X_{n,k}\distr\pi_{\eta,\,g_{\Delta_p}}\qquad
  \text{as \ $n\to\infty$.}
 \]
\end{Thm}

For the proof of Theorem \ref{THM:PD}, we use the following lemma.

\begin{Lem}\label{LEM:PD}
Let \ $\{\eta_n:n\in\ZZ_+\}$ \ be extended real--valued measures on
 \ $\Delta_p$ \ such that \ $\eta_n(\Delta_p\setminus\Lambda_r)<\infty$ \ for
 all \ $n,r\in\ZZ_+$. \ Then the following statements are equivalent:
 \begin{enumerate}
 \item[(a)] \ $\eta_n(x+\Lambda_r)\to\eta_0(x+\Lambda_r)$ \ as \ $n\to\infty$
        \ for all \ $r\in\NN$, \ $x\in\Delta_p\setminus\Lambda_r$,
 \item[(b)] \ $\DS\int_{\Delta_p}f\,\dd\eta_n\to\int_{\Delta_p}f\,\dd\eta_0$
        \ as \ $n\to\infty$ \ for all \ $f\in\cC_0(\Delta_p)$.
 \end{enumerate}
\end{Lem}

\noindent{\bf Proof.} By Theorem \ref{Portmanteau}, (b) is
equivalent to
 \begin{enumerate}
  \item[$\text{(b}^\prime\text{)}$]
   $\eta_n|_{\Delta_p\setminus U}\weak\eta_0|_{\Delta_p\setminus U}$ \ as
    \ $n\to\infty$ \ for all \ $U\in\cN_e$ \ with \ $\eta_0(\partial U)=0$.
 \end{enumerate}
Obviously, if
 \ $\eta_n|_{\Delta_p\setminus U}\weak\eta_0|_{\Delta_p\setminus U}$
 \ holds for some \ $U\in\cN_e$ \ with \ $\eta_0(\partial U)=0$ \ then
 \ $\eta_n|_{\Delta_p\setminus V}\weak\eta_0|_{\Delta_p\setminus V}$
 \ holds for all \ $V\in\cN_e$ \ with \ $V\supset U$ \ and
 \ $\eta_0(\partial V)=0$.
\ Since \ $\{\Lambda_r:r\in\NN\}$ \ is an open neighbourhood basis
of \ $e$ \ and \ $\partial\Lambda_r=\emptyset$ \ for all
  \ $r\in\ZZ_+$, \ $\text{(b}^\prime\text{)}$ is equivalent to
 \begin{enumerate}
  \item[$\text{(b}^{\prime\prime}\text{)}$]
   $\eta_n|_{\Delta_p\setminus\Lambda_r}
    \weak\eta_0|_{\Delta_p\setminus\Lambda_r}$
   \ as \ $n\to\infty$ \ for all \ $r\in\NN$.
 \end{enumerate}
For distinct elements \ $x,y\in\Delta_p$, \ let \ $\varrho(x,y)$ \
be the number \ $2^{-m}$, \ where \ $m$ \ is the least nonnegative
integer for which \ $x_m\not=y_m$. \ For all \ $x\in\Delta_p$, \ let
\ $\varrho(x,x):=0$. \ Then \ $\varrho$ \ is an invariant metric on
\ $\Delta_p$ \ compatible with the topology of \ $\Delta_p$ \ (see
Theorem 10.5 in Hewitt and Ross \cite{HR}). Let \
$d(x,y):=\sum_{k=0}^\infty2^{-k}\bone_{\{x_k\not=y_k\}}$ \ for all
 \ $x,y\in\Delta_p$. \ Then \ $d$ \ is a metric on \ $\Delta_p$ \ equivalent to \
$\varrho$, \ since \ $\varrho(x,y)\leq d(x,y)\leq2\varrho(x,y)$ \
for all \ $x,y\in\Delta_p$. \ Hence the original topology
 of \ $\Delta_p$ \ and the topology on \ $\Delta_p$ \ induced by
 the metric \ $d$ \ coincide. Then weak convergence of bounded
 measures on the locally compact group \ $\Delta_p$ \ can be considered as weak
 convergence of bounded measures on the metric space \ $\Delta_p$ \
 equipped with the metric \ $d$.

We show that the set
  \ $M:=\{\bone_{x+\Lambda_c}:c\in\NN,\,x\in\Delta_p\}$ \
 is convergence determining for the weak convergence of probability measures on
 \ $\Delta_p$. \ For this one can check that Proposition 4.6 in Ethier and Kurtz
\cite{EK} is applicable with the following choices: \ $S:=\Delta_p$
\ equipped with the metric \ $d$, \ $S_k$ \ is the set \
$\{0,1,\ldots,p-1\}$ \ for all \ $k\in\NN,$ \
 $d_k$ \ is the discrete metric on \ $S_k,$ $k\in\NN,$ \ and
\[
  M_k:=\{f_{c_k}:c_k\in S_k\},\quad k\in\NN,
    \quad\text{where}\quad
   f_{c_k}(x):=\begin{cases}
                   1 & \text{if \ $x=c_k$,} \\
                   0 & \text{if \ $x\ne c_k$,}
                \end{cases}
   \quad x\in S_k,\;\; k\in\NN.
\]
 For checking we note that for each \ $c\in\NN$ \ and \ $x\in\Delta_p$, \
 the function \ $\bone_{x+\Lambda_c}$ \ is bounded and continuous,
 since the set \ $x+\Lambda_c$ \ is open and closed. Moreover, for
 each \ $k\in\NN$, \ $S_k$ \ with the metric \ $d_k$ \ is
 a complete separable metric space.

It is easy to check that \ $M$ \ is also a convergence determining
set for the weak convergence of bounded measures on \ $\Delta_p$. \
Consequently, $\text{(b}^{\prime\prime}\text{)}$ is equivalent to
 \begin{enumerate}
  \item[$\text{(b}^{\prime\prime\prime}\text{)}$]
   $\DS\int_{\Delta_p}
        \bone_{x+\Lambda_c}\,\dd\eta_n|_{\Delta_p\setminus\Lambda_r}
    \to\int_{\Delta_p}
        \bone_{x+\Lambda_c}\,\dd\eta_0|_{\Delta_p\setminus\Lambda_r}$
    \ as \ $n\to\infty$ \ for \ all \ $x\in\Delta_p$, \; $c,r\in\NN$.
 \end{enumerate}
Clearly, this is equivalent to
 \begin{enumerate}
  \item[$\text{(b}^{\prime\prime\prime\prime}\text{)}$]
   $\DS\eta_n\big((x+\Lambda_c)\cap(\Delta_p\setminus\Lambda_r)\big)
    \to\eta_0\big((x+\Lambda_c)\cap(\Delta_p\setminus\Lambda_r)\big)$
    \ as \ $n\to\infty$ \ for all \ $x\in\Delta_p$ \ and for all
    \ $c,r\in\NN$.
 \end{enumerate}
We have
 \[
  (x+\Lambda_c)\cap(\Delta_p\setminus\Lambda_r)
  =\begin{cases}
    \Lambda_c\setminus\Lambda_r
     & \text{if \ $r\geq c$ \ and \ $x\in\Lambda_c$,} \\
    \emptyset & \text{if \ $r<c$ \ and \ $x\in\Lambda_r$,} \\
    x+\Lambda_c & \text{otherwise.}
   \end{cases}
 \]
If \ $r\geq c$ \ then \ $\Lambda_c\setminus\Lambda_r$ \ can be
written as a union of \ $p^{r-c}-1$ \ disjoint sets of the form
  \ $y+\Lambda_r$ \ with \ $y\in\Lambda_c\setminus\Lambda_r$. \
  Consequently, $\text{(b}^{\prime\prime\prime\prime}\text{)}$ and (a) are
equivalent. \proofend

\noindent{\bf Proof.} \emph{(Proof of Theorem \ref{THM:PD})} \ The
local mean of any random element with values in \ $\Delta_p$ \ is \
$e$
 \ (with respect to the local inner product \ $g_{\Delta_p}=0$\,).
Moreover, for each \ $U\in\cN_e$, \ there exists \ $r\in\ZZ_+$ \
such that
 \ $\Lambda_r\subset U$.
\ Hence, in view of Theorem \ref{THM:GAISER}, it is enough to check
that
 \begin{enumerate}
 \item[$\text{(i}^\prime\text{)}$]
  \ $\DS\max_{1\sleq k\sleq K_n}\PP(X_{n,k}\in\Delta_p\setminus\Lambda_r)\to0$
   \ as \ $n\to\infty$ \ for all \ $r\in\ZZ_+$,
 \item[$\text{(ii}^\prime\text{)}$]
  \ $\DS\sum_{k=1}^{K_n}\EE\,f(X_{n,k})\to\int_{\Delta_p}f\,\dd\eta$ \ as
   \ $n\to\infty$ \ for all \ $f\in\cC_0(\Delta_p)$.
 \end{enumerate}
Clearly
 \ $\{x\in\Delta_p:(x_0,x_1,\dots,x_d)\not=0\}=\Delta_p\setminus\Lambda_{d+1}$,
 \ hence $\text{(i}^\prime\text{)}$ and (i) are identical.
Applying Lemma \ref{LEM:PD} for \
$\eta_n:=\sum_{k=1}^{K_n}\PP_{X_{n,k}}$
 \ and \ $\eta_0:=\eta$, \ we conclude that $\text{(ii}^{\prime\prime}\text{)}$
 and (ii) are equivalent.
\proofend

\begin{Rem}
Theorem \ref{THM:GenRademacher} has the following consequence on \
$\Delta_p$. \ If \ $x_n\in\Delta_p$, \ $n\in\NN$ \ such that \
$x_n\to e$, \ and
 \ $\{X_{n,k}:n\in\NN,\,k=1,\ldots,K_n\}$ \ is a rowwise i.i.d.\ array of
 random elements in \ $\Delta_p$ \ such that \ $K_n\to\infty$ \ and
 \ $\PP(X_{n,k}=x_n)=\PP(X_{n,k}=-x_n)=\frac{1}{2}$, \ then the array
 \ $\{X_{n,k}:n\in\NN,\,k=1,\ldots,K_n\}$ \ is infinitesimal and
 \ $\sum_{k=1}^{K_n}X_{n,k}\distr\delta_e$.
\end{Rem}

\section{Limit theorems on the \ $p$--adic solenoid}
\label{gencltS}

Let \ $p$ \ be a prime. The \ $p$--adic solenoid is a subgroup of \
$\TT^\infty$, \ namely,
 \[
  S_p=\big\{(y_0,y_1,\dots)\in\TT^\infty:
            \text{$y_j=y_{j+1}^p$ \ for all \ $j\in\ZZ_+$}\big\}.
 \]
This is a compact Abelian \ $T_0$--topological group having a
countable basis
 of its topology.
Its character group is
 \ $\hS_p=\{\chi_{d,\ell}:d\in\ZZ_+,\,\ell\in\ZZ\}$, \ where
 \[
  \chi_{d,\ell}(y):=y_d^\ell,\qquad
  y\in S_p,\quad d\in\ZZ_+,\quad\ell\in\ZZ.
 \]
The set of all quadratic forms on \ $\hS_p$ \ is
 \ $\qq_+\big(\hS_p\big)=\{\psi_b:b\in\RR_+\}$, \ where
 \[
  \psi_b(\chi_{d,\ell}):=\frac{b\ell^2}{p^{2d}},\qquad
  d\in\ZZ_+,\quad\ell\in\ZZ,\quad b\in\RR_+.
 \]
An extended real--valued measure \ $\eta$ \ on $S_p$ is a
 L\'evy measure if and only if $\eta(\{e\})=0$ and
 \ $\int_{S_p}(\arg y_0)^2\,\dd\eta(y)<\infty$.
\ The function \ $g_{S_p}:S_p\times\hS_p\to\RR$,
 \[
  g_{S_p}(y,\chi_{d,\ell}):=\frac{\ell h(\arg y_0)}{p^d},\qquad
  y\in S_p,\quad d\in\ZZ_+,\quad\ell\in\ZZ,
 \]
 is a local inner product for \ $S_p$.

Theorem \ref{THM:GAISER} has the following consequence on the
 \ $p$--adic solenoid \ $S_p$.

\begin{Thm}[Gauss--Poisson limit theorem]\label{THM:GPS}
Let \ $\{X_{n,k}:n\in\NN,\,k=1,\ldots,K_n\}$ \ be a rowwise
independent array
 of random elements in \ $S_p$.
\ Suppose that there exists a quadruplet \
$(\{e\},a,\psi_b,\eta)\in\cP(S_p)$
 \ such that
 \begin{enumerate}
 \item[(i)] \ $\DS\max_{1\sleq k\sleq K_n}
          \PP(\exists\,j\leq d:|\arg((X_{n,k})_j)|>\vare)\to0$
        \ as \ $n\to\infty$ \ for all \ $d\in\ZZ_+$, \ $\vare>0$,
 \item[(ii)] \ $\DS\exp\left\{\frac{i}{p^d}
                        \sum_{k=1}^{K_n}\EE\,h(\arg((X_{n,k})_0))\right\}
           \to a_d$
        \ as \ $n\to\infty$ \ for all \ $d\in\ZZ_+$,
 \item[(iii)] \ $\DS\sum_{k=1}^{K_n}\Var\,h(\arg((X_{n,k})_0))
          \to b+\int_{S_p}h(\arg(y_0))^2\,\dd\eta(y)$ \ as \ $n\to\infty$,
 \item[(iv)] \ $\DS\sum_{k=1}^{K_n}\EE\,f(X_{n,k})\to\int_{S_p}f\,\dd\eta$
        \ as \ $n\to\infty$ \ for all \ $f\in\cC_0(S_p)$.
 \end{enumerate}
Then the array \ $\{X_{n,k}:n\in\NN,\,k=1,\ldots,K_n\}$ \ is
infinitesimal and
 \[
  \sum_{k=1}^{K_n}X_{n,k}\distr\delta_a*\gamma_{\psi_b}*\pi_{\eta,\,g_{S_p}}
  \qquad\text{as \ $n\to\infty$.}
 \]
\end{Thm}

If the limit measure has no generalized Poisson factor \
$\pi_{\eta,\,g_{S_p}}$
 \ then the truncation function \ $h$ \ can be omitted.
The proof can be carried out as in case of Theorem \ref{THM:GT}.

\begin{Thm}[CLT]\label{THM:GS}
Let \ $\{X_{n,k}:n\in\NN,\,k=1,\ldots,K_n\}$ \ be a rowwise
independent array
 of random elements in \ $S_p$.
\ Suppose that there exist an element \ $a\in S_p$ \ and a
nonnegative real
 number \ $b$ \ such that
 \begin{enumerate}
 \item[(i)] \ $\DS\exp\left\{\frac{i}{p^d}
                      \sum_{k=1}^{K_n}\EE\,\arg((X_{n,k})_0)\right\}
           \to a_d$
        \ as \ $n\to\infty$ \ for all \ $d\in\ZZ_+$,
 \item[(ii)] \ $\DS\sum_{k=1}^{K_n}\Var\,\arg((X_{n,k})_0)\to b$ \ as
        \ $n\to\infty$,
 \item[(iii)]
     \ $\DS\sum_{k=1}^{K_n}\PP(\exists\,j\leq d:|\arg((X_{n,k})_j)|>\vare)\to0$
        \ as \ $n\to\infty$ \ for all \ $d\in\ZZ_+$, \ $\vare>0$.
 \end{enumerate}
Then the array \ $\{X_{n,k}:n\in\NN,\,k=1,\ldots,K_n\}$ \ is
infinitesimal and
 \[
  \sum_{k=1}^{K_n}X_{n,k}\distr\delta_a*\gamma_{\psi_b}.
 \]
\end{Thm}

Theorem \ref{THM:GenRademacher} has the following consequence on \
$S_p$.

\begin{Thm}[Limit theorem for rowwise i.i.d.\ Rademacher array]
\label{THM:TRademacherS} Let $x^{(n)}\in S_p$, \ $n\in\NN$ \ such
that $x^{(n)}\to e$.
 Let $\{X_{n,k}:n\in\NN,\,k=1,\ldots,K_n\}$ be a rowwise i.i.d.\ array
of random elements in \ $S_p$ \ such that \ $K_n\to\infty$ \ and
 \[
  \PP(X_{n,k}=x^{(n)})=\PP(X_{n,k}=-x^{(n)})=\frac{1}{2}.
 \]
Then the array \ $\{X_{n,k}:n\in\NN,\,k=1,\ldots,K_n\}$ \ is
infinitesimal.

If \ $b$ \ is a nonnegative real number then
 \[
  \sum_{k=1}^{K_n}X_{n,k}\distr\gamma_{\psi_b}
  \qquad\Longleftrightarrow\qquad
  K_n\big(\arg(x_0^{(n)})\big)^2\to b.
 \]
\indent Moreover,
 \[
  \sum_{k=1}^{K_n}X_{n,k}\distr\omega_{S_p}
  \qquad\Longleftrightarrow\qquad
  K_n\big(\arg(x_0^{(n)})\big)^2\to\infty.
 \]
\end{Thm}

\noindent {\bf Acknowledgements.} This research has been carried out
during the stay of the authors at the Mathematisches
For\-schungs\-ins\-ti\-tut Oberwolfach supported by the ``Research
in Pairs'' Programme. The first author has been supported by the
Hungarian Scientific Research Fund
 under Grant No.\ OTKA--F046061/2004. The second author has been
 supported by the Polish Scientific Research Fund, Grant 1PO3A o31 29.
 The first and third authors have been supported by the Hungarian
 Scientific Research Fund under Grant No.\ OTKA--T048544/2005.
\proofend

\end{document}